\documentclass[12pt]{amsart}
\usepackage{amssymb, amscd, amsmath, amsthm, latexsym, enumerate}

\renewcommand{\geq}{\geqslant}
\renewcommand{\leq}{\leqslant}

\newtheorem{theorem}{Theorem}
\newtheorem{lemma}[theorem]{Lemma}
\newtheorem{corollary}[theorem]{Corollary}

\theoremstyle{definition}

\newtheorem{maintheorem}{Theorem}

\newcommand{\bY}{\breve Y}

\newcommand{\F}{\mathbb F}

\newcommand{\Z}{\mathbb Z}

\DeclareMathOperator{\cd}{cd}
\DeclareMathOperator{\gd}{gd}
\DeclareMathOperator{\Wh}{Wh}

\DeclareMathOperator{\cone}{cone}
\DeclareMathOperator{\Der}{Der}

\DeclareMathOperator{\Ext}{Ext}
\DeclareMathOperator{\Hom}{Hom}

\DeclareMathOperator{\inc}{inc}
\DeclareMathOperator{\lk}{lk}
\DeclareMathOperator{\PR}{Pr}

\newcommand{\ol}[1]{{\overline{#1}}}

\newcommand{\wt}[1]{{\widetilde{#1}}}

\newtheorem*{cor*}{Corollary}

\begin{document}

\title{aspherical $PD_4$-pairs}

\author{James F.  Davis \and J. A. Hillman }
\address{{Department of Mathematics, Indiana University,}
\newline
{Bloomington, IN 47405 USA} 
\newline
{School of Mathematics and Statistics, University of Sydney,}
\newline
Sydney,  NSW 2006, Australia }

\email{jfdavis@indiana.edu, jonathan.hillman@sydney.edu.au}

\begin{abstract}
We consider which groups $\pi$ and $PD_3$-complexes $Y$ are realized by 
$PD_4$-pairs $(X,Y)$ with $X$ aspherical and $\pi_1X\cong\pi$, 
and show that such a pair may be assembled from $(D^4,S^3)$ 
and $PD_4$-pairs of groups by adding 1-handles and mapping cylinders 
of $\Z\pi_1N$-homology equivalences over boundary components $N$ if and only if $H^2(\pi;\Z\pi)=0$.
(This includes all such pairs with $\pi_1$-injective boundary components
and all with $\pi$ a free group,  but none with $\cd\pi=2$.  
Adding a 1-handle includes connected sum.)
If $\pi$ has a finite 2-dimensional $K(\pi,1)$ complex it is realizable by some pair $(X,Y)$,
but there are no obvious building blocks analogous to $PD_4$-pairs of groups,
except for when $\pi$ is a $PD_2$-group.
\end{abstract}


\maketitle

In \cite{dh2} we show that a compact aspherical 4-manifold $M$ 
with elementary amenable fundamental group is determined up to homeomorphism
by $\partial{M}$ and the map $\pi_1\partial M \to \pi_1M$,
and that either $\pi = \pi_1M$ is trivial, $\Z$, a Baumslag-Solitar group $BS(1,m)$, a polycyclic $PD_3$-group, or a polycyclic $PD_4$-group.
The significance of the restriction ``elementary amenable" is that this is the largest class of groups 
(of finite cohomological dimension) for which surgery techniques and the $s$-cobordism theorem are known.
(The successful implementation of these techniques in \cite{dh2} is based on \cite{dh1}.)

A closed 3-manifold $N$ is the boundary of a compact aspherical 4-manifold if and only if 
it has no 2-sided projective plane and its only proper summands with finite fundamental group 
are copies of the Poincar\'e  homology sphere $S^3/I^*$.
(See Theorem \ref{asph4mfdbdry} below.)
Thus the key questions are:
``which groups $\pi$ are fundamental groups of compact aspherical 4-manifolds?"
and ``given such a group $\pi$, 
what are the possible boundaries for 4-manifolds $M\simeq{K(\pi,1)}$?".
To move beyond the class of elementary amenable groups (at present)
we must work on the level of homotopy equivalence, rather than homeomorphism,
and so  it is natural to consider $PD_4$-pairs as well as 4-manifolds.

We shall assume throughout that $(X,\partial{X})$ is a $PD_4$-pair 
such that $X$ is aspherical, and shall write $\pi=\pi_1X$.
(Some of our results may be stated in terms of manifolds, where it seems appropriate.)
We assume also that $\pi$ is non-trivial and $\partial{X}$ is non-empty,
as the case $\pi=1$ is easy and the case $\partial{X}=\emptyset$
reduces to the study of $PD_4$-groups (which is a question of algebra).
This paper is organized around the cohomological dimension of $\pi$
which is at most 3, since the boundary is nonempty.

A pair consisting of a group $\pi$ and a $PD_3$-complex $Y$ 
is realizable by a $PD_4$-pair $(X,\partial{X})$ 
with $X\simeq{K(\pi,1)}$ and $\partial{X}=Y$
if and only if there is a map $p:Y\to{K(\pi,1)}$ 
which satisfies some obviously necessary conditions, 
detailed in Theorem \ref{PD4pair}.
We show also that two such pairs are homotopy equivalent if and only if their peripheral systems
(enhanced by orientation data) are equivalent, in Theorem \ref{C-PD4}.

The exemplary case is that of $PD_4$-pairs of groups $(\pi,\mathcal{K})$, 
where $\mathcal{K}$ is a finite set of monomorphisms into $\pi$.
This essentially algebraic notion corresponds to $PD_4$-pairs $(X,\partial{X})$ with $X$ 
and its boundary components aspherical, 
and the inclusions of the boundary components being $\pi_1$-injective.
If the latter condition holds it is easily seen that $H^2(\pi;\Z\pi)=0$.
We shall show that if $H^2(\pi;\Z\pi)=0$ then $(X,\partial{X})$ 
may be constructed from copies of $(D^4,S^3)$ and $PD_4$-pairs of groups,
by the addition of 1-handles and mapping cylinders $MCyl(f_i)$
of $\Z\pi_1N_i$-homology equivalences $f:Y_i\to{N_i}$ over boundary components $N_i$.

The condition $H^2(\pi;\Z\pi)=0$ holds if $\pi$ is free, 
and more generally if it is a connected sum 
of duality groups of dimension 3 and a free group.
However, when $\cd\pi=2$ the inclusion of the boundary is never $\pi_1$-injective,
and when $H^2(\pi;\Z\pi)$ is not cyclic there are no obvious ``standard" examples.

If $\cd\pi\leq2$ then every $PD_4$-pair $(X,\partial{X})$ 
with $X\simeq{K(\pi,1)}$ has connected boundary and
is homotopy equivalent to a boundary connected sum of  $PD_4$-pairs
$(X_i,N_i)$ with $X_i$ aspherical and $N_i$ an indecomposable $PD_3$-complex 
(Theorems \ref{pi=F(r)} and \ref{cd2asph}). 
More precisely, we have the following results.

\begin{maintheorem} \label{A}
If $\pi\cong{F(r)}$ then $\partial{X}$ is connected,
and $(X,\partial{X})\cong(\natural_{i=1}^k(E(r_i)\cup{MCyl(f_i))}\natural(C,\Sigma)$,
where $E(r_i)\cong\natural^{r_i}D^3\times{S^1}$ or $\natural^{r_i}D^3\tilde\times{S^1}$,
$f_i:N_i\to\partial{E(r_i)}$ is a $\Z{F(r_i)}$-homology equivalence,
and $N_i$ is indecomposable,  for $i\leq{k}$,  
$\sum_{i\leq{k}}r_i=r$
and $C$ is the cone on a $PD_3$-homology $3$-sphere $\Sigma$.
\end{maintheorem}

In Section \ref{free_group}, we give examples of indecomposable 3-manifolds which are $\Z{F(r)}$-homology equivalent to
proper connected sums $\#^r(S^2\times{S^1})$, with $r>1$.

\begin{maintheorem} \label{B}
If $\cd\pi=2$  then $\partial{X}$ is connected, and 
$(X,\partial{X})\simeq(\natural_{i=1}^k(X_i,N_i))\natural(C,\Sigma)$, 
where $(X_i,N_i)$ is a $PD_4$-pair such that $X_i$ is aspherical
and $N_i$ is either aspherical, $S^2\times{S^1}$ or $S^2\tilde\times{S^1}$, 
for $i\leq{k}$, 
and $C$ is the cone on a $PD_3$-homology $3$-sphere $\Sigma$.
\end{maintheorem}

This result has a stronger formulation when $H^2(\pi;\Z\pi)\cong\Z$,
and so $\pi$ is a $PD_2$-group (see Theorem \ref{PD2cor}).
Every group $\pi$ with a finite 2-dimensional $K(\pi,1)$ complex is the fundamental
group of a compact aspherical 4-manifold, 
but there is an unresolved gap between the condition $\cd\pi\leq2$
and having geometric dimension 2, 
reflected in the Eilenberg-Ganea Conjecture.

Our third main result considers groups $\pi$ such that
$\cd\pi=3$ and $H^2(\pi;\Z\pi)=0$,
and summarizes the essential points of
Theorems \ref{PD4duality}, \ref{*3dg} and \ref{APD factors}.

\begin{maintheorem} \label{C}
If $\cd\pi=3$ and $H^2(\pi;\Z\pi)=0$ then $\pi\cong(*_{j=1}^rG_j)*F(s)$,
where $G_j$ is a duality group of dimension $3$,  for all $j\leq{r}$,  and $r>0$.
Assume that $\partial{X}$ has $m$ components $Y_i$,
and let $\nu_i$ be the image  of $\pi_1Y_i$ in $\pi$, for $i\leq{m}$.
Let $\mathcal{K}$ be the set of inclusions of the $\nu_i$ in $\pi$.
Then
\begin{enumerate}
\item{}if $r=1$ and $s=0$ then $(\pi,\mathcal{K})$ is a $PD_4$-pair of groups;
\item{} at most $r+s-1$ of the subgroups $\nu_i$ are not $PD_3$-groups;
\item{} $(X,\partial{X})$ is homotopy equivalent to a pair obtained by adding mapping cylinders of 
$\Z\nu_i$-homology equivalences over boundary components $Y_i'$
of an aspherical $PD_4$-pair $(X',Y')$ with $\pi_1$-injective boundary;
\item{}If  $\pi_1\inc_i$ is injective for $i\leq{m}$
then $(X,\partial{X})$ may be built up from $(D^4,S^3)$ and $PD_4$-pairs 
by boundary connected sums and attaching $1$-handles.
\end{enumerate}
\end{maintheorem}

If $\pi$ is a duality group, as in (3), then the components of $Y'$ are aspherical,
but in general they are connected sums of aspherical $PD_3$-complexes.

We have only partial results for more general groups $\pi$ with $\cd\pi=3$.
There are examples with $\pi$ indecomposable and  $H^2(\pi;\Z\pi)\not=0$,
and there are also duality groups of dimension 3 (such as $BS(1,m)\times\Z$)
which have no subgroups which are $PD_3$-groups, 
and so cannot be the fundamental group of an aspherical $PD_4$-pair.

The first section presents our notation and terminology for homological group theory,
Poincar\'e duality and connected sums.
In \S2 we discuss the notion of peripheral system and give the basic consequences of 
Poincar\'e duality for aspherical $PD_n$-pairs.
We then give a brief outline of the lower dimensional cases, 
and of the decomposition of $PD_3$-complexes, 
which we shall use in our analysis of boundary components.
In \S4 we give a criterion for a $PD_4$-pair to be aspherical, 
and find restrictions on the $PD_3$-complexes which may be boundary components of such pairs.
In \S5 we  give a technical result allowing for the decomposition of aspherical pairs into simpler pieces.
Theorem 10,  the first of the results mentioned above, appears in \S6, 
and we show that the enhanced peripheral system 
is a complete invariant for the homotopy type of an aspherical $PD_4$-pair.
The next two sections consider the cases $\cd\pi=1$ and $\cd\pi=2$.
Theorem \ref{A} is proven in  \S7 (see Theorem \ref{pi=F(r)}),  and Theorem \ref{B} is proven in \S8
(see Theorem \ref{cd2asph}).
The next four sections give our results on the case $\cd\pi=3$.
Part (1) of Theorem \ref{C} is in \S9 (see Theorem \ref{PD4duality}),
parts (2) and (3) are in \S10 (see Theorem \ref{*3dg}), 
and part (3) in \S11 (see Theorem \ref{APD factors}).
The final brief section has  a short list of questions that remain open.


\section{notation and terminology}

A group $G$ is infinite if and only if $H^0(G;\Z{G})=0$, 
and has one end if and only if $H^i(G;\Z{G})=0$ for $i=0$ and 1 (see \cite{DD89}).  
If $G$ is nontrivial, finitely generated and torsion free and is not a proper free product 
then either $G$ has one end or $G\cong\Z$.
The {\em cohomological dimension} of $G$ is 
\[
\cd G = \sup~\{q \mid H^q(G; A) \not = 0 \text{ for some $\Z{G}$-module $A$}. \}
\]
A group $G$ {\em has type $FP$} if the augmentation left $\Z G$-module $\Z$ has 
a finite length resolution by finitely generated projective $\Z G$-modules, 
and {\em is of type $F$} it is the fundamental group of a finite aspherical complex.
If $G$ is of type $FP$ then the Euler characteristic  
$\chi(G)=\sum_{i\geq0}(-1)^i\dim H_i(G;\mathbb{Q})$ is well-defined.
We shall say that $\gd{G}=n$ if $\cd{G}=n$ and there is a finite $n$-dimensional $K(G,1)$ complex.
If $n\not=2$ then $\gd{G}=n\Leftrightarrow\cd{G}=n$,
but it remains an open question whether every finitely presented group
$G$ with $\cd G= 2$ has $\gd{G}=2$.

The group $G$ is a {\em duality group} of dimension $n$
if it has type $F$ (it is more standard to assume that $G$ has type $FP$, but this nuance is not relevant for our paper and we simplify the exposition) and $H^q(G; \Z G)=0$  for $q \not = n$.
A duality group of dimension $n$ is torsion-free and has cohomological dimension $n$.  
Since $\Z{G}$ is a $(\Z{G},\Z{G})$-bimodule, 
$H^n(G; \Z{G})$ is a right $\Z{G}$-module.   
The right $\Z{G}$-module  $H^n(G; \Z G)$ is called {\em the dualizing module} for $G$.

A duality group $\pi$ is a $PD_n$-group (a {\em Poincar\'e duality group} 
of dimension $n$) if  $H^n(\pi; \Z\pi) =\Z$ (as an abelian group).
The corresponding homomorphism $w: \pi \to \Z^\times$ is called
the {\it orientation character\/} of the Poincar\'e duality group $\pi$.
Let $\Z^w$ be the $(\Z\pi,\Z\pi)$-bimodule which is additively the infinite cyclic group 
but with  $g . n = n.g =  w(g) n$,
and let $\varepsilon_w:\Z\pi\to\Z^w$ be the additive extension of $w$ 
(the {\it $w$-twisted augmentation homomorphism}).
Then $H^n(\pi; \Z\pi) =\Z^w$ as right $\Z\pi$-modules.

A $PD_n$-group $\pi$ satisfies Poincar\'e duality,
in that $H^q(\pi;\Z\pi)$ is isomorphic to $H_{n-q}(\pi;\Z\pi)$ as an abelian group,
for all $q\geq0$.
However, if we wish an isomorphism of (left) $\Z\pi$-modules, 
we need to replace $H^q(\pi;\Z\pi)$ by the conjugate module,
defined in terms of $w$.
If $B$ is a right $\Z\pi$-module, let $\ol{B}$ be the left module with the same underlying group 
and $\Z\pi$-module structure determined by $g.b=w(g)rg^{-1}$, for all $b\in{B}$ and $g\in\pi$. 
If $L$ is a left $\Z\pi$-module then $\Hom_{\Z\pi}(L,\Z\pi)$ is a right module, and so the {\it conjugate dual}
$L^\dagger=\ol{\Hom_{\Z\pi}(L,\Z\pi)}$ is again a left module.
If $M$ is a finitely generated free left $\Z\pi$-module then $M^\dagger$ is free, 
and is {\it non}-canonically isomorphic to $M$.
There are Poincar\'e duality isomorphisms
$\ol{H^q(\pi;\Z\pi)}\cong H_{n-q}(\pi;\Z\pi)$ (see \cite{dh1}).
If $\rho$ is a subgroup of finite index in a torsion-free group $\pi$,
then $\rho$ is a $PD_n$-group if and only if $\pi$ is.  

Poincar\'e duality groups were first studied by Johnson and Wall, 
and duality groups are due to Bieri and Eckmann; 
a basic reference is the book of Brown \cite{Brown}.

In this paper a $PD_n$-complex is a Poincar\'e complex of dimension $n$ 
which has the homotopy type of a finite CW-complex and a $PD_n$-pair is a Poincar\'e pair 
of dimension $n$ which has the homotopy type of a finite CW-pair.  
An {\it aspherical $PD_n$-pair\/} (or {\it $PD_n$-group with boundary})
is a $PD_n$-pair $(X,Y)$ with $X$ aspherical.
If the components of $Y=\sqcup_{i=1}^mY_i$ are also aspherical and $\pi_1$-injective,
then $(X,Y)$ is the geometric realization of a {\it $PD_n$-pair of groups},
and we may refer to such a pair as a ``$PD_4$-pair of groups", for brevity.
(A subspace $V\subset{U}$  is {\it $\pi_1$-injective\/} if the inclusions
of the components of $V$ induce monomorphisms on fundamental groups.)
We shall also use the algebraic notation $(\pi,\mathcal{K})$,
for $PD_n$-pairs of groups, as in \cite{BE}.
(Here $\mathcal{K}$ is the set of inclusions $\pi_1Y_i\to\pi=\pi_1X$, for $i\leq{m}$.)

The notion of ``$PD_n$-group with boundary" is between ``$PD_n$-pair"
and ``$PD_n$-pair of groups".
(See  \cite{Gi22} for an algebraic formulation of the case when $\pi$ 
is a duality group of dimension $n-1$ and the boundary components are aspherical, 
but not necessarily $\pi_1$-injective.)

We next review the notions of connected sum,  attaching 1-handles, and boundary connected sum and attaching 1-handles, generalizing all of these notions to the Poincar\'e category.   The fundamental constructions are decompositions and triads.  An {\em $n$-manifold decomposition} $(X;X_0,X_1)$ is an $n$-manifold $X$ (without boundary) expressed as a union $X = X_0 \cup X_1$ where $X_0$ and $X_1$ are manifolds with boundaries $\partial X_0 = X_0 \cap X_1 = \partial X_1$.  An {\em $n$-manifold triad} $(X;Y_0,Y_1)$ is an $n$-manifold $X$ with subsets $Y_0, Y_1 \subset \partial X$ so that  $(\partial X; Y_0,Y_1)$ is a decomposition.    If $(X;Y_0,Y_1)$ and $(X';Y'_0,Y'_1)$ are $n$-manifold triads and $h : Y_1 \xrightarrow{\cong} Y_0'$ is a homeomorphism, then $(X \cup_h X'; Y_0, Y_1')$ is an $n$-manifold triad.   An $n$-manifold decomposition is obtained by gluing  triads $(X_0; \emptyset, \partial X_0)$ and $(X_1; \partial X_1, \emptyset)$.
If the two  triads are oriented and $h$ is orientation-reversing, then the glued triad is oriented.    

There are parallel definitions of a Poincar\'e complex decomposition and a Poincar\'e triad, as well as gluing Poincar\'e triads, but here one can glue along a homotopy equivalence, replacing $X \cup_h X'$ by a double mapping cylinder if desired.   (In the definition of a Poincar\'e triad $(X; Y_0,Y_1)$, it is a matter of taste where $Y_i$ is a subcomplex of $X$ or if there are specified maps $Y_i \to X$.)

For any  $n$-manifold $X$ (without boundary), there is a decomposition $(X;\ol X, D)$ where $D$ is an $n$-disc.   The disc theorem says that this is unique up to homeomorphism when $X$ is connected.   For a connected Poincar\'e complex $X$ of dimension $n$, there a homotopy equivalence $X \simeq \ol X \cup D$  and a decomposition $(\ol X \cup D; \ol X, D)$ where $D$ is an $n$-disc and $\ol X$ has the homotopy type of a complex of dimension $n - 1$ if $n > 3$ and $X$ is finitely dominated by a 2-complex if $n = 3$.   The disc theorem \cite[Theorem 2.4]{Wall_P} says that for $X$ connected, this decomposition exists and is unique up to homotopy type.

A connected sum $X \# X'$ of $n$-manifolds (without boundary) is obtained by gluing the triads $(\ol X;  \emptyset, \partial D)$ and  $(\ol X'; \partial D',\emptyset)$ along a homeomorphism $h : \partial D \xrightarrow{\cong} \partial D'$.  A connected sum $X \# X'$ of Poincar\'e complexes is obtained by the same formula, but gluing along a homotopy equivalence.  (See also \cite{Wall_P}.)

Let $(X,Y)$ be an $n$-manifold with boundary.   Attaching a 1-handle $X \cup h^1$ is the result of gluing manifold triads
$$
(D^1 \times D^{n-1}; D^1 \times S^{n-2}, S^0 \times D^{n-1}) \cup_h (X;h(S^0 \times D^{n-1}),Z) 
$$
where $h : S^0 \times D^{n-1} \hookrightarrow Y$ is an embedding so that $Z = \ol{Y - h(S^0 \times D^{n-1})}$ is a manifold with boundary.  There is a parallel notion of attaching a 1-handle to a Poincar\'e pair $(X,Y)$  along an embedding
$h : S^0 \times D^{n-1} \hookrightarrow \sqcup D_i$ where $Y_i$ are the components of $Y$ and $Y_i \simeq \ol{Y_i} \cup D_i$  are Wall's decompositions.

There are three different flavors of attaching a 1-handle; either the image of $h$ lies in different  components of $X$, or in the same component of $X$ but in different components of $Y$, or in different components of $Y$.  The first flavor is the {\em boundary connected sum}.  Here $(X'', Y'') = (X,Y) \natural (X', Y')$ where $Y'' = (Y\setminus{Y_i})\sqcup(Y'\setminus{Y_j'})\sqcup(Y_i\#Y_j')$ and $X'' \simeq X \vee X'$.   The second flavor is the {\em self-boundary connected sum} where if $(X', Y') = (X,Y) \cup h^1$, then $X' \simeq X \vee S^1$ and $Y$ is a connected sum of components of $Y$.   The last flavor is a boundary connected sum of $(X,Y)$ with either $D^{n-1} \times S^1$ or $D^{n-1} \tilde\times S^1$.   In all cases if $(X,Y)$ is aspherical, then so is $(X,Y) \cup h^1$. 

We note finally that such connected sums generally depend on choices of components and local orientations for the discs. 
However we have chosen not to complicate the notation here.

\section{General remarks on the peripheral systems of aspherical $PD_n$-pairs}

Basic invariants of a  $PD_n$-pair include its fundamental group, 
the number of boundary components, 
and the fundamental groups of the boundary components.     
We formalize this by introducing the notion of a peripheral system 
and then study this applied to aspherical $PD_n$-pairs.

Let $(X,\partial{X})$ be a $PD_n$-pair with boundary
$\partial{X}=\sqcup_{i=1}^m\partial_i{X}$.  
We will usually assume that the ``ambient space" $X$ is connected, 
except when considering assembling such pairs by forming connected sums of several such pairs.
However, we will not assume the boundary is connected.
Let $\inc_i : \partial_i{X} \to{X}$ be the inclusion of the $i$-th boundary component.   
(If  $\partial{X}$ is connected we shall drop the label $i=1$.) 
Choosing base points $x_0 \in{X}$ and $x_i\in \partial_jX$
and paths $\gamma_i$ from $x_i$ to $x_0$, 
define the {\em fundamental group system of $X$} to be the tuple 
$(\pi_1X ,  \pi_1\inc_1, \dots , \pi_1\inc_m)$.   
The {\em orientation character} of $X$ is the homomorphism 
$w =w_1X  : \pi_1X\to \Z^{\times}= \{\pm 1\}$.  
The {\em peripheral system of $X$} is the fundamental group system together with the orientation character.

An isomorphism 
$(G,  j_1: G_1 \to G, \dots, j_m : G_m \to G, w) \to (G', j'_1: G'_1 \to G', \dots, j'_m : G'_m \to G',w') $
of peripheral systems consists of a permutation $\sigma \in S_k$, 
group isomorphisms $\theta : G \to G'$, $\theta_i : G_i \to G'_{\sigma(i)}$, 
and elements $g_i \in G$ so that $w= w' \circ \theta$ and 
$\theta_i \circ j_i \circ c_{g_i} = j'_i \circ \theta_{\sigma(i)}$ where $c_{g_i}$ is conjugation by $g_i$.   
The isomorphism class of the peripheral system is an invariant of $X$, 
independent of the choice of base points and paths.

If $X$ is aspherical  then $\pi=\pi_1X$ is of type $F$, 
and if $\partial{X}=\emptyset$ then $\pi$ is a $PD_n$-group with orientation character $w=w_1X$.
If $\partial{X}$ is non-empty then $\pi$ is no longer a $PD_n$-group, 
and $w_1X$ may not be determined by $\pi$.
In this case there are Poincar\'e-Lefshetz duality isomorphisms
$\overline{H^q(X;\Z\pi)}\cong{H_{n-q}(X,\partial{X};\Z\pi)}$,
where the conjugation uses the anti-involution induced by $w=w_1X$.

\begin{lemma} 
\label{components of boundary}
Let $(X,\partial{X})$ be an aspherical $PD_n$-pair and let $\pi=\pi_1X$. Then
\label{boundary components}
\begin{enumerate}
\item \label{cdpin} $\cd \pi = n \Longleftrightarrow \partial{X}$ is empty.
\item \label{cdpin-1}$\cd \pi <n\Longleftrightarrow \partial{X}$ has at most 
$\beta_{n-1}(\pi;\mathbb{F}_2)+1$ components.
\item \label{pdn-1} $ \pi$ is $PD_{n-1} \Longrightarrow \partial{X}$ has one or two components.
\end{enumerate}
\end{lemma}

\begin{proof}
\eqref{cdpin}
If $\partial{X}=\emptyset$ then 
$H^n(\pi;\mathbb{F}_2)=H^n(X;\mathbb{F}_2)= H_0(X;\F_2) = \F_2 \not = 0$.  
If $\partial{X}$ is non-empty, then 
\[
H^n(\pi;\Z\pi)=H^n(X;\Z\pi)=H_0(X,\partial{X};\Z\pi)=H_0(\wt{X}, \partial \wt{X}; \Z) = 0, 
\]
 and so $\cd \pi<n$.

\eqref{cdpin-1}   and \eqref{pdn-1} follow from the exact sequence 
\[
H_1(X,\partial{X};\mathbb{F}_2)\to
{H_0(\partial{X};\mathbb{F}_2)}\to
{H_0(X;\mathbb{F}_2)} \to H_0(X,\partial{X}; \F_2), 
\]
which is isomorphic to 
\[
H^{n-1} (\pi;\F_2) \to H_0(\partial{X}; \F_2) \to \F_2 \to H^n(\pi; \F_2),
\]
by Poincar\'e-Lefschetz duality with $\F_2$ coefficients.
\end{proof}

An example to keep in mind is $M^l \times D^{n-l}$, where $M$ is a closed aspherical manifold.      
Another example is a compact contractible manifold with boundary a homology sphere.   
Note that the boundary components of an aspherical manifold need not be aspherical.   
The bound $\beta_{n-1}(\pi;\mathbb{F}_2)+1$ on the number of boundary components
in part (2) can be achieved when $\pi$ is a free product of $PD_{n-1}$-groups,
for example when $M$ is a compact surface with many open disks removed.
Note also that the punctured torus $T_o$ and the thrice-punctured sphere $D_{oo}$ 
each have fundamental group $F(2)$, 
and so the boundary is not determined by $\pi$.

Now we move from the discussion of the number of boundary components to 
their fundamental groups. 
We shall write $H_i(X)$ and $H_i(G)$ for the integral homology of a space $X$ or group $G$, for simplicity.

\begin{theorem}
\label{cd<n}
Let $(X,\partial{X})$ be an aspherical $PD_n$-pair.
Then $\pi=\pi_1X$ is of type $F$ and $\cd\pi\leq{n}$.
In particular, $\pi$ is torsion free.
\begin{enumerate}
\item{}If $q<{n-1}$ then $H^q(\inc):H^q(\pi;\Z\pi)\to{H^q(\partial{X};\Z\pi)}$ is an isomorphism,
and there is an exact sequence 
\[0\to{H^{n-1}(\pi;\Z\pi)}\to{H^{n-1}(\partial{X};\Z\pi)}\to\Z\to0.
\]
\item If $\cd \pi < n$ and $\pi \not = 1$, then the image of $\pi_1\inc_i$ is infinite for all $i$.
\item If $\cd \pi \leq n-2$, then $\partial{X}$ is connected and $\pi_1\inc$ is an epimorphism.
\item$\ker(\pi_1\inc)$ is a perfect group if and only if $H^{n-2}(\pi;\Z\pi)=0$. 
\item If $\cd \pi = n-2$, then $\pi_1\inc$ has nontrivial kernel.
\item If $\pi$ is a $PD_{n-2}$-group, then $\ker(\pi_1\inc)$ 
has abelianization  $\Z$.
\item If $\pi$ is a duality group of dimension $n-1$ then $\ker(\pi_1\inc_i)$ is acyclic and
$\mathrm{Im}(\pi_1\inc_i)$ is a $PD_{n-1}$-group, for each boundary component $\partial_iX$.
\item If $\pi$ is a $PD_{n-1}$-group with orientation character  $w_{\pi}$,
then each $\ker(\pi_1\inc_i)$ is a perfect group and either
\begin{enumerate}
\item $\partial{X}$ has two components and each $\pi_1\inc_i$ is an epimorphism and
\[
w_1X = w_\pi : \pi \to \Z^{\times}
\] 
or
\item $\partial{X}$ has one component, 
$[\pi:\mathrm{Im}(\pi_1\inc)]=2$ and
\[
w_1X = w_{\pi} \cdot w_{\partial} : \pi \to \Z^{\times}
\]
where $w_\partial: \pi \to \Z^{\times}$ has kernel $\mathrm{Im}(\pi_1\inc)$.
\end{enumerate}
\end{enumerate}
\end{theorem}

\begin{proof}
The first assertion follows from the long exact sequence of cohomology with coefficients $\Z\pi$
for the pair $(X,\partial{X})$, Poincar\'e-Lefschetz duality 
\[\ol{H^q(X,\partial{X}; \Z\pi)} \cong{H_{n-q}(X; \Z\pi)},
\] 
asphericity of $X$, and the fact that $H_0(\pi;\Z\pi)\cong\Z$ and $H_q(\pi;\Z\pi)=0$ if $q>0$.
Similarly, 
\begin{equation} \label{PLD}
 \ol{H^q(\pi; \Z\pi)} =\ol{H^q(X; \Z\pi)} \cong \wt H_{n-q-1}(\partial {X};\Z\pi))=
 \wt H_{n-q-1}(\partial \wt{X}),
\end{equation}
for all $q$.

(2) Take $j = 0$.      
Note that since $X=K(\pi,1)$ is finite-dimensional, $\pi$ is torsion-free.  
If $\pi \not = 1$, then $\pi$ is infinite so $H^0(\pi;\Z\pi) = (\Z\pi)^{\pi} = 0$.   
Thus $H_{n-1}(\partial \wt{X})= 0$.  
Hence $\partial \wt{X}$ has no closed components, 
as we can see by either replacing $\Z$ by $\F_2$ or passing to the orientation double cover of $X$.

(3) Take $q = n-1$.   
If $\cd \pi < n-1$, then $\wt H_0(\partial \wt{X}) = 0$, and hence $\pi_1\inc$ is an epimorphism.

(4),(5),(6)   Taking $q = n-2$,
we see $H_1(\partial \wt{X}) \cong \ol{H^{n-2}(X; \Z\pi)}= \ol{H^{n-2}(\pi; \Z\pi)}$.  
 If $\cd \pi= n-2$, this is nonzero and if $\pi$ is a  $PD_{n-2}$-group, then $H_1(\partial \wt{X}) = \Z$.  
If $\cd \pi \leq n-3$, then $H_1(\partial \wt{X}) = 0$.  
In all cases, $H_1(\partial \wt{X})$ is the abelianization of $\ker(\pi_1\inc)= \pi_1\partial \wt{X}$.

(7)  If $\pi$ is a duality group of dimension $n-1$ then 
\[
H_q(X, \partial{X}; \Z\pi)=H^{n-q}(X;\Z\pi)=0\quad{\mathrm{for~all}}\quad{ q>1}.
\]
Hence $H_q(\partial{X};\Z\pi)=0$ for all $q>0$, and so $\ker(\pi_1\inc_i)$ 
is acyclic, for each boundary component  $\partial_iX$.

If $N$ is a component of the preimage of $\partial_iX$ in $\widetilde{X}$ then 
$N$ is acyclic, and the projection of $N$ onto $\partial_iX$ is a regular covering,  
with covering group $\nu$ the image of $\pi_1N$ in $\pi$.
Hence the cellular chain complex $C_*(\partial_iX;\Z\nu)$
is a finitely generated projective resolution for $\Z$.
Since $H^{n-1}(\partial_iX;\Z\nu)\cong{H_0(\partial_iX;\Z\nu)=\Z}$,
by Poincar\'e duality for $\partial_iX$, it follows that $\nu$ is a $PD_{n-1}$-group.

(8)  Take $j = n-1, n-2$.   Taking $j = n-1$ we see $\wt H_0(\partial \wt{X}) = \Z$, 
and thus $\partial \wt{X}$ has two components. 
Thus either $\partial{X}$ has two components with 
$\pi_1\inc_i$ an epimorphism or  $\partial{X}$ has one component and
the image of $\pi_1\inc$ has index 2.   
In either case, taking $j = n-2$, we see that $H_1(\partial \wt{X}) = 0$, 
so that the each map $\pi_1\inc_i$ has a perfect kernel.   
In either case, there is an exact sequence of $\Z\pi$-modules
$$
0 \to H_1(X,\partial{X}; \Z\pi) \to H_0(\partial{X}; \Z \pi) \xrightarrow{\beta} H_0(X; \Z\pi) \to 0.
$$

If $\pi$ is a $PD_{n-1}$-group and $\partial{X}$ has two components, 
then 
\[
\ol{H^{n-1}(X;\Z\pi)} \cong H_1(X,\partial{X}; \Z\pi)
 \cong \ker(\beta) \cong \Z.
 \] 
Thus  $\ol{\Z^{w_\pi}}\cong\Z$,  so $w_{\pi_1X} = w_\pi$ as claimed.

Suppose finally that $\pi$ is a $PD_{n-1}$-group and that $\partial{X}$ 
has one component with 
$\rho := \mathrm{Im}(\pi_1\inc)$ having index 2 in $\pi$.   
Then  
\[
\ol{H^{n-1}(X;\Z\pi)} \cong H_1(X,\partial{X}; \Z\pi) \cong \ker(\beta)\cong \Z^{w_\partial}.
\] 
Thus $\ol{\Z^{w_\pi}}\cong\Z^{w_\partial}$, 
so  $w_1X = w_\pi \cdot w_\partial$ as claimed.
\end{proof}

In particular, if  $\cd \pi\leq{n-2}$ or if $\pi$ is a $PD_{n-1}$-group and $\partial{X}$ has two components, then $\partial{X}$ is orientable if and only if $X$ is orientable.  

If $n\geq4$ there are compact $n$-manifolds representing each type considered in the theorem.
In case (8)(b) we may take a closed aspherical $(n-1)$-manifold $N$ 
with a double cover $\widehat{N} \to{N}$ and set 
$M = \widehat{N} \times_{\Z/2} [-1,1]$.   
Note that in this case $M\simeq{N}$, 
but $M$ and $N$ must have different orientation characters.   

We do not have a result comparable to (8) for more general groups 
of cohomological dimension $n-1$, 
but see Theorem \ref{*3dg} below for the case when $n=4$ and $\pi$ is a proper free product of 
duality groups of dimension 3.

The following is a consequence of Lemma \ref{components of boundary}
and Theorem \ref{cd<n}.

\begin{corollary}
\label{pi-inc inj}
Let $(X,\partial{X})$ be an aspherical $PD_n$-pair with non-empty boundary,
and let $\pi=\pi_1X$.
If  $\cd\pi\leq{n-2}$ (so $\partial{X}$ is connected) and $\pi_1\inc$ is a monomorphism
then $\cd\pi\leq{n-3}$ and $\pi_1\inc$ is an isomorphism.
If $\cd\pi=n-1$ and $\pi_1\inc_i$ is a monomorphism for all $i\leq{m}$ 
then $H^{n-2}(\pi;\Z\pi)=0$.
\qed
\end{corollary}

In particular, if $\cd\pi=n-2$ there are no examples with $\pi_1$-injective boundaries. 
This is already clear when $n=2$ and $\pi=1$!

\section{lower dimensions}

By way of warm-up, we outline the situation when the dimension  $n$ is 2 or 3.

Every $PD_2$-pair other than $S^2$, $\mathbb{RP}^2$ or $(D^2,S^1)$ is  a $PD_2$-pair of groups,
since it is represented by a compact aspherical 2-manifold with $\pi_1$-injective aspherical boundary.
  
If $(X,\partial{X})$ is an aspherical $PD_3$-pair and $\pi=\pi_1X$ is finite 
then $\pi=1$ and $\partial{X}=S^2$;
if $\pi$ is infinite then asphericity of the ambient space and Poincar\'e duality imply that 
$H_1(\partial{X};\Z\pi)\cong{H^1(\pi;\Z\pi)}$ and $H_2(\partial{X};\Z\pi)=0$, 
so the components of the boundary are aspherical.
Moreover $\chi(Y)=2\chi(\pi)$, so $\chi(\pi)\leq0$.

The following are equivalent:
\begin{enumerate}
\item$\pi$ has one end;
\item{} the inclusion of $\partial{X}$ into $K(\pi,1)$ is $\pi_1$-injective on each component;
\item $(X,\partial{X})$ is a $PD_3$-pair of groups.
\end{enumerate}
If $\pi$ has two ends then $\pi\cong\Z$ and $\partial{X}=T$ or $Kb$.
It is not known whether every $FP$ group $\pi$ with one end and such that 
$\cd\pi=2$ and $\chi(\pi)\leq0$
is the ambient group of a $PD_3$-pair of groups $(\pi,\mathcal{G})$.

An aspherical $PD_3$-pair $(X,\partial{X})$ with $\pi=\pi_1X\not=1$ is isomorphic to 
a boundary connected sum of $PD_3$-pairs of groups 
{\it up to connected sums with copies of $D^2$-bundles over $S^1$}
if and only if the boundary components are aspherical \cite[Corollary 13.11.2]{Hi20}.
If $X$ is a 3-manifold and $G$ is a proper free product we may use the Loop Theorem 
to show that $X$ is either a proper boundary connected sum or has a 1-handle, 
and so we may reduce to the case when $\pi$ is indecomposable and has one end.
(The resulting manifold pair is then an aspherical $PD_3$-pair.)
However the Algebraic Loop Theorem \cite[Chapter 4]{Hi20}
is not strong enough to establish this reduction for all $PD_3$-pairs.
It is not known whether every $PD_3$-pair is homotopy equivalent to a 3-manifold with boundary.

We shall use the fact that $PD_3$-complexes have essentially unique factorizations
into indecomposables \cite{Tu89}.
The indecomposables are analyzed in detail in \cite[Chapters 2--7]{Hi20}.
The key facts are
\begin{enumerate}
\item{} the homotopy type of a $PD_3$-complex $P$ is determined by the
{\it fundamental triple\/} 
$(\pi_1P, w_1P, c_{P*}[P])$, where $c_P:P\to{K(\pi_1P,1)}$ is the classifying map
and $[P]\in{H_3(P;\Z^{w_1P})}$ is a fundamental class \cite[Corollary 2.3.1]{Hi20};
\item{}if $P$ is a $PD_3$-complex such that $\pi_1P$ is a proper free product then $P$ 
has a corresponding decomposition as a connected sum
\cite[Theorem 2.7]{Hi20}; and
\item{}if $\pi_1P$ has no element $g$ with $g^2=1$ and $w_1P(g)=-1$ 
then $P$ is a connected sum of aspherical $PD_3$-complexes,
copies of $S^2\times{S^1}$ and $S^2\tilde\times{S^1}$,
and indecomposable $PD_3$-complexes whose fundamental groups 
are generated by elements of finite order \cite[Theorem 4.8 and Corollary 7.10.1]{Hi20}.
\end{enumerate}

\section{asphericity and boundary}

It is easy to see that a $PD_4$-complex $X$ is aspherical if and only if $\pi_0X = 1$, 
$\pi=\pi_1X$ has one end, and $\pi_2X=0$.
In the bounded case the condition on ends must be modified.
The boundary connected sum of aspherical manifolds with boundary
is again aspherical, but if neither summand is contractible then
the fundamental group of the sum is a nontrivial free product,
and so has more than one end.

In the next theorem we refer to the Bass Conjectures \cite{Ba76},
which we outline very briefly.
If $P$ is a finitely generated projective $\Z\pi$-module then
it is the image of an idempotent $n\times{n}$-matrix $A$ with entries in $\Z\pi$,
for some $n\geq0$.
The Kaplansky rank $\kappa(P)$ is the coefficient of 1 in the trace of $A$. 
The weak Bass Conjecture is the assertion that $\kappa(P)$ equals  the rank of $ \Z \otimes_{\Z\pi}P$.
An analytic argument originally due to Kaplansky (unpublished) shows
that if $\kappa(P)=0$, then $P=0$.
See  \cite{BV98} for a proof of this positivity result, 
in terms of traces of idempotent square matrices with coefficients
in the von Neumann algebra $\mathcal{N}\pi$.

\begin{theorem}
\label{asphcrit}
Let  $(X,\partial{X})$ be a $PD_4$-pair with boundary $\partial{X}=\sqcup_{i=1}^m\partial_iX$
and  such that $\pi = \pi_1X\not=1$.
Then $X$ is aspherical if and only if  $\pi_1\partial_iX$ 
has infinite image in $\pi$ for each $i\leq{m}$, 
the inclusions induce an isomorphism 
$H^1(\pi;\Z\pi)\cong \oplus_{i=1}^m {H^1(\pi_1\partial_iX; \Z\pi)}$, 
and $\pi_2X=0$.
If  $\cd\pi\le2$ then $X$ is aspherical if and only if
$\partial{X}\not=\emptyset$,
$(X,\partial{X})$ is $1$-connected, and $\chi(X)=\chi(\pi)$.
\end{theorem}

\begin{proof}
If $X$ is aspherical then $H_*(X;\Z\pi)=H_* (\wt X) = 0$ for $* \not = 0$.
In particular, $\pi_2X=H_2(X; \Z\pi)  = 0$.  
Poincar\'e-Lefschetz duality and the long exact cohomology sequence of a pair imply 
$H^q(X;\Z\pi) \cong \oplus_{i=1}^mH^q(\partial_iX; \Z\pi)$ for $q= 0, 1, 2$.  
Using $q = 1$, we have an isomorphism
$H^1(\pi;\Z\pi)\to\oplus_{i=1}^mH^1(\pi_1\partial_iX;\Z\pi)$. 
Since $X$ is aspherical and $\pi\not=1$, 
the group $\pi$ is infinite, and so $H^0(X;\Z\pi)=0$.
Hence $H^0(\partial_iX;\Z\pi)=0$, 
and so the image of $\pi_1\partial_iX$ in $\pi$ is infinite for $i\leq{m}$.

Conversely,  if $\pi$ is infinite then $H_4(X;\Z\pi)=0$, 
while the first two conditions imply that
$H_3(X;\Z\pi)=H^1(X,\partial{X};\Z\pi) = 0$.
Hence if also $\pi_2X=0$ then $X$ is aspherical.

Suppose now that $\cd\pi\leq2$.
If $X$ is aspherical then $\chi(X)=\chi(\pi)$, 
$\partial{X}\not=\emptyset$ (by Lemma \ref{components of boundary}),
and $(X,\partial{X})$ is $1$-connected
(by Theorem \ref{cd<n}).

If these conditions hold then $H^3(X;\Z\pi)=H_1(X,\partial{X};\Z\pi)=0$, 
by Poincar\'e duality and the fact that
$(X,\partial{X})$ is $1$-connected.
Since $\partial{X}\not=\emptyset$, 
the cellular chain complex $C_*(X;\Z\pi)$
is chain homotopy equivalent to a finite free chain complex $C_*$
of length $3$.   
Let $C^* = \Hom_{\Z\pi}(C_*,\Z\pi)$ be the associated cochain complex. 
Since $H^3(X;\Z\pi)=0$, $C^2 \to C^3$ is a split surjection, 
hence $C_3 \to C_2$ is a split injection.   
Thus $H_3(X; \Z\pi) = 0$ and $C_*$ is chain homotopy equivalent to 
a finite free chain complex $D_*$ of length $2$.

Let $Z_q\leq D_q$ be the submodule of $q$-cycles.
The generalized Schanuel's Lemma (see e.g. Lemma VIII.4.4 of \cite{Brown}) 
and the fact that $\cd \pi \leq 2$ imply that $Z_1$ is projective.    
Since $H_1 (D_*) = H_1(\wt X) = 0$, 
it follows that $D_2\cong{Z_1}\oplus{Z_2}$. 
Hence $D_*\cong{E_*}\oplus \Pi$, 
where $E_*$ is a projective resolution of $\mathbb{Z}$ and $\Pi=Z_2$ 
is a projective module concentrated in degree 2.
On tensoring with $\mathbb{Z}$ and using the condition $\chi(X)=\chi(\pi)$,
we see that $\mathbb{Z}\otimes_{\Z\pi}\Pi=0$.
Since $\cd\pi\leq2$, the Bass conjectures hold for $\pi$ \cite{Ec86},
and so $\Pi=0$.
Hence $X$ is aspherical.
\end{proof}

If $X$ has $\pi_1$-injective boundary then 
$H^2(\pi;\mathbb{Z}\pi)=H_2(X,\partial{X};\Z\pi)=0$.
Retracts of groups of type $FP$ are again of type $FP$ \cite[Proposition 3.4]{Le15}.
Hence the indecomposable free factors of $\pi$ are either infinite cyclic or
are duality groups of dimension 3.
In particular, if $\pi$ has finitely many ends then either $X\simeq{S^1}$ or 
$\pi$ is a duality group of dimension 3 and the components of $\partial{X}$ are aspherical.
In the latter case the homotopy type of $X$ is determined by its peripheral system, 
which is a $PD_4$-pair of groups.

A {\it $PD_3$-homology sphere\/} is a $PD_3$-complex with perfect fundamental group.
It is an easy consequence of the argument for Theorem \ref{cd<n} that
if $(X,\partial{X})$ is a $PD_4$-pair then the number of indecomposable summands 
of components of $\partial{X}$ which are not $PD_3$-homology spheres is bounded 
in terms of $\pi$ alone.
However, we may increase the number of homology sphere summands without 
changing $\pi$ by taking boundary connected sums with suitable contractible 4-manifolds.

\begin{theorem}
\label{asph4mfdbdry}
Let $Y$ be a $PD_3$-complex.
If $Y$ is a component of $\partial{X}$ for some $PD_4$-pair with $X$ aspherical 
then the indecomposable summands of $Y$ are either aspherical, or are
copies of $S^2\times{S^1}$ or $S^2\tilde\times{S^1}$, or are $PD_3$-homology spheres.
\end{theorem}

\begin{proof}
Let $\pi=\pi_1X$ and $w=w_1X$.  

The inclusion $\inc:Y\to{X}$ induces an exact sequence
\begin{equation*}
1 \to \kappa \to \pi_1Y \xrightarrow{\pi_1(\inc)} \pi
\end{equation*}
where $\pi := \pi_1X$ and $\kappa := \ker (\pi_1(\inc))$.
We claim that $H_1\kappa = \kappa^{ab}$ is $\Z$-torsion-free.   Indeed
\begin{equation*}
H_1\kappa = H_1(Y;\Z\pi) \subset H_1(\partial{X}; \Z\pi) \xleftarrow{\cong} 
H_2(X,\partial{X}; \Z\pi) \cong \ol{H^2(\pi;\Z\pi)}
\end{equation*}
is torsion-free,
since the second cohomology group is torsion-free for any finitely presented group
(see  \cite[Proposition 13.7.1]{Ge}). 

Suppose that $Y = Y_1 \# Y_2$.
Thus $\pi_1Y =  \pi_1Y_1 * \pi_1Y_2$.   
If $\pi_1X_1$ is generated by elements of finite order then $\pi_1N_1 \subset \kappa$,
since $\pi$ is torsion-free.
Thus $H_1X_1 \to H_1\kappa \to H_1Y$ is the zero map, 
since the left group is finite and the middle group is torsion-free.   
But $H_1Y = H_1Y_1 \oplus H_1Y_2$.   Thus $H_1Y_1 = 0$.

The fact that $\pi$ is torsion-free implies also that if 
$g\in\pi_1Y$ has order 2 then  $w_1Y(g)=0$, since $w_1Y=w\pi_1\inc$.
Thus $Y$ is a connected sum of aspherical $PD_3$-complexes, 
copies of $S^2\times{S^1}$ or $S^2\tilde\times{S^1}$,
and $PD_3$-complexes whose fundamental groups 
are generated by elements of finite order.
(See Theorems 2.7, 4.8, 6.8, 6.10 and 7.10 of \cite{Hi20}.)
\end{proof}

The next corollary is formulated in terms of manifolds, 
as this is needed for the implication $(1)\Rightarrow(2)$.

\begin{corollary}
\label{asph4mbc}
Let $N$ be a closed $3$-manifold.
Then the following are equivalent
\begin{enumerate}
\item$N$ is a connected sum of aspherical $3$-manifolds and copies of
$S^2\times{S^1}$, $S^2\tilde\times{S^1}$ and $S^3/I^*$;
\item$N=\partial{M}$ for some aspherical $4$-manifold $M$;
\item$N$ is a component of $\partial{M}$ for some aspherical $4$-manifold $M$.
\end{enumerate}
If $M$ is a compact aspherical $4$-manifold and $\partial {M}$ has $k>1$ components
then we may attach 1-handles to get a compact aspherical $4$-manifold with group $\pi*F(k-1)$
 and connected boundary.
\end{corollary}

\begin{proof}
Every closed aspherical 3-manifold bounds a compact aspherical 4-manifold, 
by the relative hyperbolization theorem of \cite{DJW}.
The $S^2$-bundle spaces $S^2\times{S^1}$ and $S^2\tilde\times{S^1}$ bound the corresponding $D^3$-bundle spaces ${D^3\times{S^1}}$ and $D^3\tilde\times{S^1}$,
while the homology 3-sphere $S^3/I^*$ bounds a contractible 4-manifold \cite{FQ}. 
These 4-manifolds are clearly also aspherical.
The boundary connected sum $M=M_1\natural\dots\natural{M_n}$ 
of a finite set of aspherical 4-manifolds
$M_1,\dots,M_n$ with connected boundaries is again aspherical, and
$\partial{M}=\partial{M_1}\#\dots\#\partial{M_n}$. Thus (1) implies (2).
Clearly (2) implies (3),
while indecomposable homology 3-spheres are either aspherical or homeomorphic to $S^3/I^*$,
and so (3) implies (1), by Theorem \ref{asph4mfdbdry}.

The final assertion is clear.
\end{proof}

The implication $(1)\Rightarrow(3)$ does not require \cite{DJW}, 
since every 3-manifold $P$ is a boundary component of $P\times[0,1]$, 
and we may again perform iterated boundary connected sums to realize
a connected sum $N=\#_{k\in{K}}P_k$ 
as a component of the boundary of a compact aspherical 4-manifold, 
but now one with disconnected boundary.

\section{decomposing $PD_n$-pairs into simpler pieces}

One of our goals is to show that aspherical $PD_n$-pairs can be assembled from simpler pieces.
In order to do this,
we shall need to extend the notion of boundary connected sum.
If $(U,V)$ and $(U',V')$ are finite CW-pairs such that $U$ and $V$ are connected and the subspaces $V$ and $V'$ are  
$PD_{n-1}$-complexes, for some $n>1$, 
then we can form a ``boundary connected sum" $(U,V) \natural (U',V') := (U \vee {U'},V\#V')$.
Similarly, we may attach a 1-handle to  such pair $(U,V)$.
Of course we must choose components and local orientations for the discs involved,
but we do not wish to complicate the notation by recording such choices.

Theorem C of \cite{dh1} states that a finite CW-pair $(U,V)$ with $U$ aspherical and $V$ a Poincar\'e complex of dimension $n-1$ is a Poincar\'e pair if and only if 
$H^i(U,V; \Z\pi)$ vanishes for $i \not = n$ and is an infinite cyclic abelian group if $i = n$.   Here $\pi = \pi_1U$.

The next lemma provides a criterion for recognizing when $(U,V)$ is a $PD_n$-pair of 
groups.

\begin{lemma}
\label{PD4paircriteria}
Let $(U,V)$ be a finite, $\pi_1$-injective CW-pair with $U$ aspherical and $V$  nonempty.  
Then $(U,V)$ is a $PD_{n}$-pair of groups if and only if 
\begin{enumerate}
\item $V$ is a disjoint union of aspherical  $PD_{n-1}$-complexes; and 
\item $\pi$ is a duality group of dimension $n-1$; and 
\item $H^{n-1}(U,V;\Z\pi)=0$ and  $H^n(U,V;\Z\pi)\cong\Z$.
\end{enumerate}
\end{lemma}

\begin{proof}
Let $(U,V)$ be a finite, $\pi_1$-injective CW-pair with $U$ aspherical and $V$ a nonempty $PD_{n-1}$-complex.  Then $H_i(U;\Z\pi) = 0 = H_i(V; \Z\pi)$ for $i \not = 0$ and $H_i(U,V;\Z\pi) = 0$ for $i \not =1 $.

If $(U,V)$ is a $PD_n$-pair of groups with $V$ nonempty, then (1) holds by definition and (2), (3) hold by Poincar\'e-Lefschetz duality.

Conversely if (1) and (2) hold, the exact sequence 
$$
\dots   \to  H^{i-1}(V; \Z\pi) \to H^i(U,V; \Z\pi) \to  H^i(U; \Z\pi) \to \dots
$$
shows that $H^i(U,V; \Z\pi) = 0$ for $i \not = n-1, n$. Then  (3) and Theorem C of \cite{dh1} show that $(U,V)$ is a Poincar\'e pair.
\end{proof}

A $\pi_1$-injective CW-pair $(U,V)$ is an  {\em aspherical CW-pair with $PD_{n-1}$ boundary},
or simply an {\em $APD_{n-1}\partial$-pair}, if $U$ and the components of $V$ are aspherical and $V$ is $PD_{n-1}$.  

\begin{theorem}
\label{PD4natural}
Let $(X,Y)=(U,V)\natural(U',V')$, where $(U,V)$ and $(U',V')$ are $APD_{n-1}\partial$-pairs,
for some $n\geq2$.
If $(X,Y)$ is a $PD_n$-pair, then so are each of $(U,V)$ and $(U',V')$.

A similar result holds if $(X,Y)$ is the result of attaching a 1-handle to
an $APD_{n-1}\partial$-pair $(U,V)$.
\end{theorem}

\begin{proof}
We  claim that if $A$ and $B$ are CW-complexes so that $F = A \vee B$ has the homotopy type of a finite CW-complex, then so does $A$ (and hence $B$).  Since $\pi_1 B$ is a retract of the finitely generated group $\pi_1F$, the group $\pi_1B$ is also finitely generated. After attaching a finite number of 2-cells to $B$ to kill its fundamental group,  $F' = A \vee B'$ has the homotopy type of a finite CW-complex with  $B'$ simply-connected.   But $B'$ is  a finitely dominated space, hence has the homotopy type of a finite complex $B''$ by \cite{Wall_FI}.  Then $F'' = A \vee \cone (B'')$ has the homotopy type of $A$ as well as the homotopy type of a finite CW-complex.  

Thus we may assume that $(U,V)$ and $(U',V')$ are finite $APD_{n-1}\partial$-pairs.

Let $\pi = \pi_1 X = \pi_1 U * \pi_1 U'$.  Consider the triple 
$$
Y = V \# V' \to \breve{Y} := V \vee V' \to X = U \vee U'
$$
which is homotopy equivalent to the triple 
$$
Y  \hookrightarrow \breve{Y}_D := (V \# V') \cup_{S^{n-2}} D^{n-1} \hookrightarrow X_D := (U \sqcup U') \cup_{\partial D^1\times  D^{n-1} } D^1 \times D^{n-1}.
$$
where the spheres and discs arise from Wall's disc theorem applied to $V$ and $V'$ and are geometrically motivated by the model of the boundary connected sum coming from adding a 1-handle to $U \sqcup U'$.  

We wish to analyze the long exact cohomology sequence of the triple $(X_D,\bY_D,Y)$, equivalently $Y \to \bY \to X$, with coefficients in a $\Z\pi$-module $M$.  Given a map $A \to X$, we use the letter $M$ to also indicate the $\Z\pi_1A$-module given by restriction.

Excision shows 
\[
H^*(\bY_D,Y ;M) \xrightarrow{\cong} H^*(D^{n-1},S^{n-2} ;M) 
\]
and 
\[
H^*(X, \bY; M) \cong H^*(U, V; M) \oplus H^*(U', V'; M).
\]
Poincare-Lefschetz duality shows 
\[
H^*(X,Y; M) \cong H_{n-*}(X; \ol M).
\]
Thus there is a long exact sequence
\begin{multline*}
 \dots \to H^{q-1}(D^{n-1},S^{n-2};M) \to H^q(U, V; M) \oplus H^q(U', V';M) \to \\ H_{n-q}(X;\ol M) \to 
H^{q}(D^{n-1},S^{n-2};M) \to \dots
\end{multline*}
It follows that for $q \not = n$, $H^q(U,V;\Z\pi) = 0$.   
 
 Since $\Z\pi$ is a free $\Z\pi_1U$-module,
\[
H^*(U, V; \Z\pi) = \prod_{|\pi : \pi_1 U|} H^*(U, V; \Z\pi_1U);
\]
similarly for $(U',V')$.    Thus for $q \not = n$, $H^q(U, V; \Z\pi_1U)$ and $H^q(U', V'; \Z\pi_1U')$ vanish.  

Next we set $n = q$ in the exact sequence above and obtain the exact sequence 
\begin{multline}\label{M}
H^{n-1}(D^{n-1},S^{n-1}; M)  \to H^n(U, V; M) \oplus H^n(U', V';M)  \\ \to H^n(X,Y; M) \to 0
\end{multline}

We first analyze this exact sequence for $M = \Z^w$.  Assume the connected sum joins the components $V_1$ of $V$ and $V_1'$ of $V'$.  Then one has the commutative square
\[
\begin{CD}
H^{n-1}(V_1; \Z^w) \oplus H^{n-1}(V_1'; \Z^w) @>>> H^{n-1}(V_1 \# V_1'; \Z^w) \\
@V\cong VV @V\cong VV \\
H^{n}(U,V; \Z^w) \oplus H^{n-1}(U,V; \Z^w) @>>> H^{n}(X,Y; \Z^w)
\end{CD}
\]
where the top horizontal arrow is an isomorphism on each summand and thus the bottom arrow is an isomorphism on each summand.  We conclude that  $H^{n}(U,V; \Z^w)$ and $H^{n}(U',V'; \Z^w)$ are both infinite cyclic.   We also see that $H^{n-1}(D^{n-1},S^{n-1}; \Z^w)  \xrightarrow{\cong} H^n(U, V; \Z^w)$ is an isomorphism, similarly for $(U',V')$.   (Note: This fact about connected sums follows from Wall's disc theorem.) 

Finally, we show  that $H^n(U,V;\Z\pi_1U)\cong\Z$.
The group $\pi$ retracts onto its subgroup $\pi_1U$ and so
acts on $\Z\pi_1U$ by right multiplication.
Let $\mathcal{L}=\overline{\Z\pi_1U}$ be the conjugate left $\Z\pi$-module.
Then $\mathcal{L}=\Z^w\otimes_\Z\Z\pi_1U$, with the diagonal left $\pi$-action, where $w : \pi \to \{\pm1\}$ is the orientation character.
If we restrict the action to $\pi_1U$ (considered as a subgroup), 
then $\mathcal{L}\cong\ol{\Z\pi_1U}$, while if we restrict it to $\pi_1U'$ then ${\mathcal{L}}\cong \sum_{\pi_1U} \Z^w$.

We have the commutative diagram 
$$
\begin{CD}
H^{n-1}(D^{n-1},S^{n-1};\mathcal{L}) @>>> H^{n-1}(U',V';\mathcal{L}) \\
@| @| \\
\prod_{\pi_1U} H^{n-1}(D^{n-1},S^{n-1};\Z^w) @>>> \prod_{\pi_1U} H^{n-1}(U',V';\Z^w)
\end{CD}.
$$
The bottom horizontal arrow is an isomorphism by the computation above, so the top horizontal arrow is also.   A diagram chase (or one of Noether's isomorphism theorems) applied to \eqref{M}  then shows that
$$
H^n(U,V; \ol{\Z\pi_1U}) \xrightarrow{\cong} H^n(X,Y; \ol{\Z\pi_1U}).
$$
The domain is isomorphic to $H^n(U,V; {\Z\pi_1U})$ and the codomain 
 is infinite cyclic by Poincar\'e-Lefschetz duality.

A similar argument  shows that
$H^n(U',V';\Z\pi_1U')\cong\Z$. Thus    $(U,V)$ and $(U',V')$ are $PD_n$-pairs by Theorem C of \cite{dh1} mentioned above.  

Suppose now that $(X,Y)$ is the result of attaching a 1-handle to
an $APD_{n-1}\partial$-pair $(U,V)$. 
If the discs used in forming the self-boundary connected sum lie in different components of
$V$ then the above argument goes through with very minor changes.
If they lie in the same component $N$ then $(X,Y)$ is the boundary connected sum of $(U,V)$ 
with either $(D^{n-1}\times{S^1},\partial)$ or $(D^{n-1}\tilde\times{S^1},\partial)$, 
and we may use the earlier argument.
\end{proof}

The next lemma is used to recognize aspherical $PD_n$-pairs which may be 
obtained from simpler  $PD_n$-pairs  by adding mapping cylinders over boundary components.

\begin{lemma}
\label{excismv}
Let $(U,A)$ and $(W,A\sqcup{Y})$ be pairs of cell complexes with
$U$ and $W$ are connected and let $X=U\cup_AW$.
If the inclusions of $A$ and $Y$ into $W$ are a homotopy equivalence
and a $\Z\pi_1W$-homology equivalence, respectively,
then $H_p(U,A;\Z\pi)\cong{H_p(X,W;\Z\pi)}\cong{H_p(X,Y;\Z\pi)}$,
for all $p$.
Similarly for cohomology.
\end{lemma}

\begin{proof}
The inclusion $(U,A)\subset(X,W)$ is excisive,  and so $H_p(U,A;\Z\pi)\cong{H_p(X,W;\Z\pi)}$.
This is in turn isomorphic to $H_p(X,Y;\Z\pi)$,
by the exact sequence of the triple $(X,W,Y)$, 
since $Y\to{W}$ is a $\Z\pi_1W$-homology equivalence.
and so $H_*(W,Y;\Z\pi)=0$.
\end{proof}

It follows that if $(U,V)$ is a pair with $U$ aspherical,
$A$ is a component of $V$ and $W=MCyl(f)$ is the mapping cylinder of
a $\Z\pi_1A$-homology equivalence $f:Y\to{A}$, and $X=U\cup_AW\simeq{U}$,
then $(X,Y)$ is a $PD_n$-pair if and only if $(U,V)$ is a $PD_n$-pair.

\section{existence: boundaries of aspherical $PD_4$-pairs}

Let $(U,\sqcup_{i=1}^mK_i)$ be a $PD_4$-pair of groups
and let $f_i:Y_i\to{K_i}$ be a $\Z\pi_1K_i$-homology equivalence, for each $i\leq{m}$.
Then the mapping cylinders $MCyl(f_i)$ are the ambient spaces of $PD_4$-pairs with
boundary $K_i\sqcup{Y_i}$.
Let $X=U\cup_{K_1}MCyl(f_1)\dots\cup_{K_k}MCyl(f_k)$
and let $Y=\sqcup_{i=1}^mY_i$.
Then $(X,Y)$ is a compact aspherical $PD_4$-pair which realizes the group $\pi_1X\cong\pi$ 
and has boundary components $\{Y_1,\dots,Y_m\}$.

The following theorem complements Theorem \ref{cd<n}.
The class of groups with $\cd=3$ considered here is limited, 
but includes all the groups in \S9 below and in \cite{dh2}.

\begin{theorem}
\label{PD4pair}
Let $\pi$ be a finitely presented group,  $w:\pi\to\mathbb{Z}^\times$ be a homomorphism, 
$Y=\sqcup_{i=1}^mY_i$ be a $PD_3$-complex and $p:Y\to{K(\pi,1)}$ be a map such that $w_1Y=p^*w$.
Let $X$ be the mapping cylinder of $p$ and let $\kappa=\ker(\pi_1p)$.
Suppose that one of the following conditions holds:
\begin{enumerate}
 \setcounter{enumi}{-1}
\item$\pi=1$, $Y$ is connected and $\pi_1Y$ is perfect;
\item$\cd \pi=1$, $Y$ is connected, $p_*=\pi_1p$ is an epimorphism, and $\kappa$ is perfect;
\item$\cd \pi=2$, $Y$ is connected,  $p_*$ is an epimorphism, and 
$H^q(p;\Z\pi)$ is an isomorphism for $q\leq2$;
\item$\pi$ is a duality group of dimension $3$ which is the ambient group of a $PD_4$-pair
of groups $(\pi,\mathcal{K})$, where $\mathcal{K}=\sqcup_{i=1}^mK_i$
and there are degree-$1$ maps $f_i:Y_i\to{K_i}$ such that $\pi_1f_j$ has perfect kernel and
$p|_{Y_i}=f_i$ for $i\leq{m}$.
\end{enumerate}
Then the pair $(X,Y)$ is a $PD_4$-pair with orientation character $w$.
\end{theorem}

\begin{proof}
The result is clear if $\pi=1$, $Y$ is connected and $\pi_1Y$ is perfect, 
for then $X$ is contractible and so $H^q(X,Y;\Z)\cong\widetilde{H}^{q-1}(Y;\Z)$ for all $q$.

When $\pi\not=1$ the image of $\pi_1Y$ in $\pi$ is infinite,
and so $H^0(Y;\Z\pi)=H^0(X;\Z\pi)=0$.
Hence $H^0(X,Y;\Z\pi)=0$.

If (1) holds then $H^1(\kappa;\Z\pi)=\Hom(\kappa,\Z\pi)=0$,
since $\kappa$ acts trivially on $\Z\pi$ and is perfect. 
It then follows from the 5-term exact sequence of low degree for $\pi_1Y$ 
as an extension of $\pi$ by $\kappa$ that $H^1(p;\Z\pi)$ is an isomorphism.
Since $H_p(X;\mathcal{R})=0$ for  any right $\Z\pi$-module $\mathcal{R}$ and $p>1$,
we have $H_4(X,Y;\Z^w)\cong{H_3(Y;\Z^w)}\cong\Z$.
Let $[X,Y]$ and $[Y]$ be compatible generators for these groups.
Since $H^q(X;\Z\pi)=0$ for $q>1$ and $H^2(Y;\Z\pi)\cong{H_1(Y;\Z\pi)} = H_1(\kappa)=0$, 
by Poincar\'e duality for $Y$,
we see that $H^q(X,Y;\Z\pi)=0$ for $q\leq3$,
and $H^4(X,Y;\Z\pi)\cong{H^3(Y;\Z\pi)}\cong\Z$.
Consideration of the commutative diagram
\begin{equation*}
\begin{CD}
0@>>>H^3(Y;\Z\pi)@>\cong>>H^4(X,Y;\Z\pi)@>>>0\\
@VVV @V\cap[Y]VV @V\cap[X,Y]VV @VVV\\
0@>>>H_0(Y;\Z\pi)@>\cong>>H_0(X;\Z\pi)@>>>0
\end{CD}
\end{equation*}
shows that cap product with the class $[X,Y]$ induces isomorphisms 
$H^q(X,Y;\Z\pi)\cong{H_{4-q}(X;\Z\pi)}$ for all $q$,
and so $(X,Y)$ is a $PD_4$-pair.

A similar argument applies in case (2), for we again have
$H_p(X;\Z\pi)=0$ and $H^q(X;\Z\pi)=0$ for $p,q>2$.
Hence $H_4(X,Y;\Z^w)\cong{H_3(Y;\Z^w)}\cong\Z$,
for $q\leq3$ and an isomorphism from $H^3(Y;\Z\pi)\cong{H^4(X,Y;\Z\pi)}$.
The hypotheses on $H^q(p)$ and the fact that $H^3(X;\Z\pi)=0$ 
imply that $H^q(X,Y;\Z\pi)=0$ for $q\leq 3$, and $(X,Y)$ is a $PD_4$-pair
as in case (1).

If (3) holds, then each of the maps $f_i$ is a $\Z\pi_1Y_i$-homology equivalence,
by Poincar\'e duality for $Y_i$ and $K_i$,
and so the claim follows from Lemma \ref{excismv}.
\end{proof}

Parts (0) and (1) of Theorem \ref{PD4pair} are covered by \cite[Proposition 11.6C]{FQ}.

The next lemma and its corollary are based on the fact that every homology 3-sphere bounds 
a contractible 4-manifold \cite{Fr82}, and so have been formulated in terms of 4-manifolds.

\begin{lemma}
\label{prime bdry}
Let $M$ be a compact aspherical  $4$-manifold and let $\pi=\pi_1M$.
If $\partial{M}$ has a component $\partial_iM=N\#\Sigma$,
where $\Sigma$ is a homology $3$-sphere and $\pi_1\Sigma\leq\ker(\pi_1\inc)$,
 then $(M,\partial{M})\simeq(M',N)\natural(C,\Sigma)$,
where $M'\simeq{M}$,  $C$ is  contractible and $\partial{C}\cong\Sigma$.
\end{lemma}

\begin{proof}
We may write $\partial_i{M}=N_o\cup\Sigma_o$, 
where $N_o$ and $\Sigma_o$ are the complements of open 3-discs in 
$N$ and $\Sigma$, respectively.
The homology 3-sphere $\Sigma$ bounds a contractible 4-manifold $C$ \cite{Fr82}.
Then $\partial{C}=\Sigma_o\cup{D^3}$.
Let $M'=M\cup_{\Sigma_o}C$. 
The inclusion of $M$ into $M'$ induces an isomorphism $\pi_1M\cong\pi_1M'$,
since $\pi_1\Sigma_o=\pi_1\Sigma$ has trivial image in $\pi$.
It also induces isomorphisms $H_i(M;\Z\pi)\cong H_i(M';\Z\pi)$ for all $i$,
since $H_i(\Sigma_o;\Z\pi)=0$ for $i>0$.
Hence this inclusion is a homotopy equivalence $M\simeq{M'}$.
\end{proof}

If $\Sigma$ is a $PD_3$-homology sphere and $C$ is the cone over $\Sigma$ then 
$C$ is contractible and $(C,\Sigma)$ is a $PD_4$-pair,
and so the analogue of Lemma \ref{prime bdry} for aspherical $PD_4$-pairs is easy.

\begin{corollary}
\label{boundarycomps}
{\it 
Let $M$ be a compact aspherical $4$-manifold and let $\pi=\pi_1M$.
If $\pi$ has no finitely generated perfect subgroup then there is a $4$-manifold 
$M'\simeq{M}$ such that no component of $\partial{M'}$ has
a homology $3$-sphere summand,
and a set of contractible $4$-manifolds $\{C_1,\dots,C_k\}$ such that 
$(M,\partial{M})\simeq(M',\partial{M'})\natural(C_1,\partial{C_1})\dots\natural(C_k,\partial{C_k})$.}
\end{corollary}

\begin{proof}
If $\partial_iM$ is a component of $\partial {M}$ then we may assume that 
$\partial_i{M}=N_i\#\Sigma_i$,
where $\Sigma_i$ is a homology 3-sphere and no summand of $N_i$ is a homology $3$-sphere.
The image of $\pi_1\Sigma_i$ in $\pi$ is trivial, 
and the claim follows from Lemma \ref{prime bdry}.
\end{proof}

We may assume further that the homology spheres $\partial{C_i}$ are irreducible,
and so none of the contractible 4-manifolds $C_i$ are proper boundary connected sums.

It follows from the Core Theorem \cite{Sc73}
(respectively, the Algebraic Core Theorem \cite{KK05})
that 3-manifold groups (respectively, coherent $PD_3$-groups) 
have no finitely generated perfect subgroups of infinite index
(although they may have infinitely generated,  perfect normal subgroups).

We show next that the homotopy types of aspherical $PD_4$-pairs
are determined by a mildly enhanced form of their peripheral data.
This relies upon the classification of  $PD_3$-complexes by their fundamental triples.
If $(X,\partial{X})$ is a $PD_4$-pair and $N$ is an orientable  component of $\partial{X}$,
let $[N]$ be the image in $H_3(\pi_1N;\Z)$ of a generator $[X,\partial{X}]$
for $H_4(X,\partial{X};\Z^w)$ under the composition of the connecting homomorphism
to $ H_3(\partial{X};\Z)$, projection onto the summand $H_3(N;\Z)$, 
and the homomorphism to $H_3(\pi_1N;\Z)$ induced by the classifying map.
The {\it enhanced\/} peripheral system of $(X,\partial{X})$ 
is the peripheral system together with the classes $[N]$.
Two such enhanced peripheral systems are equivalent if there is an isomorphism of peripheral systems which preserves each of the homology classes $[N]$, 
up to a simultaneous change of signs.

\begin{theorem} \label{C-PD4}
Let $(X,\partial{X})$ and $(X',\partial{X'})$ be aspherical $PD_4$-pairs 
with $\partial{X}$ and $\partial{X'}$ non-empty.
Then $(X,\partial{X})\simeq(X',\partial{X'})$ if and only if their enhanced peripheral systems are equivalent.
\end{theorem}

\begin{proof} 
An isomorphism of peripheral systems which preserves the enhancement
induces equivalences of the fundamental triples of corresponding boundary components,
and hence induces homotopy equivalences of such components 
which are compatible with the inclusions into the ambient spaces $X$ and $X'$.
Thus if the enhanced peripheral systems are equivalent then $(X,\partial{X})\simeq(X',\partial{X'})$.

The converse is clear.
\end{proof}

If $X$ and $X'$ are compact aspherical manifolds with orientable boundaries and 
the Farrell-Jones Conjectures hold for $\pi$ (e.g $\pi$ is elementary amenable) then we may strengthen this result to show that $(X,\partial X)$ and $(X',\partial X')$ 
that $\partial{X}$ and $\partial{X'}$ are homeomorphic pairs.
(Orientability of the boundaries is needed here because the Geometrization Theorem of 
Perelman and Thurston has not yet been proven for nonorientable  3-manifolds.)

The necessity of some level of enhancement is clear from the following example.
Let $C$ be a contractible 4-manifold with boundary the Poincar\'e homology 3-sphere
$P=S^3/I^*$, and fix an orientation for $C$.
Since self-homotopy equivalences of $P$ are orientation-preserving,
$P\#{P}$ is not homeomorphic to $P\#-\!P$.
Hence the two boundary connected sums $M_+=C\natural{C}$ and 
$M_-=C\natural-\!C$ are not homeomorphic, 
although their peripheral systems are both $\{I^**I^*\to1\}$.
In this case the enhanced peripheral systems are not equivalent.

\section{$\pi$ a free group} \label{free_group}

In this section we shall show that if $\pi$ is a free group then every $PD_4$-pair realizing $\pi$ may be 
assembled from copies of $D^3\times{S^1}$ and $D^3\tilde\times{S^1}$ by adding mapping cylinders 
of $\Z{F}$-homology equivalences and forming boundary-connected sums.
We also give some examples of such manifolds with $\pi$ nonabelian and 
indecomposable boundary.

\begin{lemma}
\label{G*Hontopi}
Let  $f:G*H\to\pi$ be a homomorphism with perfect kernel $\kappa$.
Then $f(G*H)\cong{f(G)*f(H)}$,
and $G\cap\kappa$ and $H\cap\kappa$ are perfect.
In particular, if $f(G)=1$ then $G$ is perfect.
If $\pi$ is torsion free and has no noncyclic free subgroups then one of $f(G)$ and $f(H)$ is trivial.
\end{lemma}

\begin{proof}
The homomorphism $f$ induces an epimorphism $\bar{f}:f(G)*f(H)\to\sigma=f(G*H)$.
The kernel $K$ of $\bar{f}$ is the quotient of $\kappa$ by the kernel of the canonical projection 
of $G*H$ onto $f(G)*f(H)$,  and so is perfect. 
On the other hand $K\cap{f(G)}=K\cap{f(H)}=1$, and so $K$ is free, by Proposition 1.4.5 of \cite{DD89}.
Therefore $K=1$ and $\bar{f}$ is an isomorphism.

The abelianizations $(G\cap\kappa)^{ab}=H_1(G;\Z\sigma)$ and 
$(H\cap\kappa)^{ab}=H_1(H;\Z\sigma)$
are direct summands  of $\kappa^{ab}=H_1(G*H;\Z\sigma)$ and so are trivial.
Hence $G\cap\kappa$ and $H\cap\kappa$ are perfect.
In particular, if $f(G)=1$ then $G=G\cap\kappa$ is perfect.

The final assertion is clear, 
since the free product of nontrivial torsion free groups contains noncyclic free subgroups.
\end{proof}

In particular, 
if $f$ is an epimorphism then $f(G)$ and $f(H)$ are each free products
of indecomposable factors of $\pi$.

\begin{lemma}
\label{homologyiso}
Let $f:P\to{Q}$ be a degree-$1$ map of connected $PD_3$-complexes 
such that $\kappa=\mathrm{Ker}(\pi_1f)$ is perfect.
Then $f$ is a $\Z\pi_1Q$-homology equivalence.
\end{lemma}

\begin{proof}
The homomorphisms $H_0(f;\Z\pi_1Q)$ and $H_1(f;\Z\pi_1Q)$ are isomorphisms,
since $f$ has degree 1 and $\kappa$ is perfect,
while $H^0(f;\Z\pi_1Q)$ is an isomorphism, 
since $P$ and $Q$ are connected and $f$ is an epimorphism.
As in  part (1) of Theorem \ref{PD4pair},
it follows from the five term exact sequence of low degree that
 $H^1(f;\Z\pi_1Q)$ is also an isomorphism, since
$\kappa$ acts trivially on $\Z\pi_1Q$, and is perfect.
The lemma now follows by Poincar\'e duality (with coefficients $\Z\pi_1Q)$)
using again the fact that $f$ has degree 1.
\end{proof}

The next result is essentially Theorem \ref{A} of the introduction.

\begin{theorem}
\label{pi=F(r)}
Let $(X,\partial{X})$ be a compact aspherical $PD_4$-pair such that $\pi=\pi_1X\cong{F(r)}$,
and suppose that $\partial{X}\cong(\#_{i=1}^kN_i)\#\Sigma$,
where $N_i$ is indecomposable and $H_1(N_i;\Z)\not=0$,
and $\Sigma$ is a $PD_3$-homology sphere.
Then  $(X,\partial{X})\cong(\natural(E(r_i)\cup{MCyl(f_i))}\natural(C,\Sigma)$,
where $E(r_i)\cong\natural^{r_i}D^3\times{S^1}$ or $\natural^{r_i}D^3\tilde\times{S^1}$,
$f_i:N_i\to\partial{E(r_i)}$ is a $\Z{F(r_i)}$-homology equivalence,  $\sum_{i\leq{k}}r_i=r$
and $C$ is the cone on $\Sigma$.
\end{theorem}

\begin{proof}
Since $\kappa$ is perfect, 
it follows from Lemma \ref{G*Hontopi} that the groups $\pi_1N_i$ map onto distinct nontrivial  factors $F(r_i)$ of $\pi$, for $i\leq{k}$,
and $\Sigma_{i=1}^kr_i=r$.
The induced epimorphisms $\theta_i:\pi_1N_i\to{F(r_i)}$ may be realized by 
degree-1 maps $f_i:N_i\to\partial{E(r_i)}$, where $w_1E(r_i)=w|_{F(r_i)}$
\cite[Corollary 3]{He77}.
Such maps are $\Z{F(r_i)}$-homology equivalences, by Lemma 
\ref{homologyiso}.

Let $X_i=E(r_i)\cup{MCyl(f_i)}$ and let $Y_i$ be the image of $N_i\subset{MCyl(f_i)}$ in $X_i$,
for $i\leq{k}$.
Then $(X_i,Y_i)$ is a $PD_4$-pair, 
and hence so is $(X',\partial{X'})=(\natural(X_i,Y_i))\natural(C,\Sigma)$.
The map from $X'$ to $X$ which is induced by the inclusions of the factors
$F(r_i)$ into $\pi$ is homotopic to a map which restricts to a homotopy equivalence
$(\#N_i)\#\Sigma\simeq\partial{X}$, and so $(X',\partial{X'})\simeq(X,\partial{X})$.
\end{proof}

Theorem  \ref{PD2cor} and Corollary \ref{3DGcor} below give parallel results
when $\pi$ is a $PD_2$-group or a duality group of dimension 3, respectively.

\begin{corollary}
\label{freeboundaryinj}
If $\pi$ is a free group and $\pi_1\inc$ is a monomorphism then $(X,\partial{X})$
is homotopy equivalent to a connected sum of copies of $D^3\times{S^1}$ and
$D^3\tilde\times{S^1}$.
\qed
\end{corollary}

When $\pi\cong\Z$ and $(X,Y)$ is orientable we do not need to appeal to \cite{He77}
in order to find a degree-1 map from $N$ to $S^2\times{S^1}$.
For $H^1(N;\Z)\cong[N,S^1]$ and $[N,S^2]$ maps onto 
$H^2(N;\Z)\cong[N;\mathbb{CP}^\infty]$.
If $f_q:N\to{S^q}$ are maps corresponding to generators of $H^2(N;\Z)$,
for $q=1$ and 2,  and $f=(f_2,f_1):N\to {S^2}\times{S^1}$,
then $H^*(f;\Z)$ is an isomorphism.

We shall now turn to the manifold case.
If $M$ is a compact aspherical 4-manifold then we may use Lemma  \ref{prime bdry} 
to strengthen Theorem \ref{pi=F(r)}  slightly:
$(M,\partial{M})\simeq(M',\partial{M'})\natural(C,\Sigma)$,
where $M'$ and $C$ are 4-manifolds, 
$C$ is contractible and $\partial{C}=\Sigma$.

The 3-manifolds $M(K)$ obtained by 0-framed surgeries on a knot $K$ 
with Alexander polynomial $\Delta_K(t)\dot=1$ provide examples with $\pi\cong\Z$,
for then $\pi_1M(K)$ is an extension of $\Z$ by a perfect normal subgroup,
and if $K$ is nontrivial then $M(K)$ is aspherical \cite{Ga87}.
There is then a 4-manifold $M$ such that $\pi_1M\cong\Z$ and $\partial{M}\cong{M(K)}$,
by the Existence Theorem of \cite{dh1}.

We may use Whitehead doubles of links to construct examples 
with $\pi$ a nonabelian free group and aspherical boundary.
The exterior of each component of the Whitehead link $Wh_1=5^2_1$ is a solid torus, 
in which the other component is a homologically trivial knot.
If $L$ is an $m$-component link then its Whitehead doubles 
$Wh(L)$ are obtained by replacing each component 
by a (possibly twisted) copy of such a knot.
(We do not need to specify the twisting here.)
Thus the exterior $X(Wh(L))$ is the union of $X(L)$ with $m$ copies of $X(Wh_1)$,
each attached along one of its two boundary tori.

An $m$-component link $L$ is {\it freely slice\/} if it bounds a set $\mathcal{D}_L$
of $m$ disjoint locally flat discs which are properly embedded in $D^4$,
with exterior $Z(\mathcal{D}_L)$ such that $\pi_1Z(\mathcal{D}_L)\cong{F(m)}$.

\begin{theorem}
\label{Wh(L)}
Let $L$ be an unsplittable $m$-component link with all pairwise linking numbers $0$.
If either $m=2$ or $L$ is a boundary link then $N_L=M(Wh(L))$ is aspherical, 
and bounds a compact aspherical $4$-manifold with fundamental group $F(m)$.
\end{theorem}

\begin{proof}
The hypothesis on $L$ implies that its exterior $X(L)$ is aspherical and 
the boundary is $\pi_1$-injective.
Let $E=X(Wh_1)\cup_{\lambda\times{S^1}}{D^2\times{S^1}}$ be the 3-manifold with boundary 
obtained from the exterior of the Whitehead link $Wh_1=5^2_1$ by filling in one component 
so that its longitude $\lambda$ bounds a disc.
Then it is not hard to see that the inclusion of $\partial{E}$ into $E$ is  $\pi_1$-injective,
and that $E$ is aspherical.
Since $N_L$  may be obtained from $X(L)$ by adjoining copies of $E$ 
along each boundary component, $N_L$ is aspherical,
by Van Kampen's Theorem and a Mayer-Vietoris argument.

The hypotheses on $L$ also imply that $Wh(L)$ is freely slice \cite{FT95b},
and so $Wh(L)$ bounds a set of slice discs $\mathcal{D}$ in $D^4$,
with exterior $Z_L=Z(\mathcal{D})$ and such that $\pi=\pi_1Z_L\cong{F(m)}$.
It is easy to see that $H_i(Z_L;\Z)=0$ for $i>1$.
The inclusion of $N_L=\partial{Z_L}$ sends meridians to meridians
and so defines an epimorphism from $\nu=\pi_1N_L$ to $\pi$.
Hence $\pi_1(Z_L,N_L)=0$.
These three properties imply that $Z_L\simeq\vee^mS^1$ (and so is aspherical)
\cite[Proposition 11.6C]{FQ}.
\end{proof}

The simplest non-trivial choice for $L$ is the Whitehead link $Wh_1$,
and $\Wh(Wh_1)$ has  a diagram with 24 crossings!

The kernel of any epimorphism $\nu\to{F(m)}$ 
which induces an isomorphism on abelianization is 
the intersection of the terms of the lower central series of $\nu$,
and so does not depend on the choice of epimorphism \cite{St65}.

It is not known whether every aspherical $PD_4$-pair $(X,N)$
with $\pi$ is a nonabelian free group is homotopy equivalent to a 4-manifold pair $(M,\partial{M})$.
However if this is so then $M$ is unique up to TOP $s$-cobordism \cite[Theorem 11.6A]{FQ}.

We conclude this section with two examples of aspherical smooth 4-manifolds 
with free fundamental group and aspherical boundary, 
and for which we do not need to invoke surgery.
They are each based on the ribbon knot of Figure 1.4 of \cite{HiA},
which has associated slice disc $\mathcal{D}$ such that 
$\pi_1Z(\mathcal{D})\cong\Z$.
This knot is in fact the Kinoshita-Terasaka knot $KT=11_{42n}$,  
and so is nontrivial.
The 4-manifold $Z(\mathcal{D})$ is aspherical, 
as in Theorem \ref{Wh(L)},
and its boundary $M(KT)$ is aspherical \cite{Ga87}.

\setlength{\unitlength}{1mm}
\begin{picture}(90,63)(5,-5)

\put(15,20){\line(0,1){6}}
\qbezier(15,20)(15,11)(24,11)
\qbezier(15,26)(15,35)(24,35)
\put(19,20){\line(0,1){6}}
\qbezier(19,20)(19,15)(24,15)
\qbezier(19,26)(19,31)(24,31)

\put(24,35){\line(1,0){10}}
\put(25,10){\line(0,1){23}}
\qbezier(25,10)(25,6.5)(28.5,6.5)
\put(25,36){\line(0,1){5}}
\qbezier(25,41)(25,50)(34,50)
\qbezier(29,41)(29,46)(34,46)
\qbezier(28.5,6.5)(32,6.5)(32,10)
\put(29,11){\line(1,0){7}}
\put(29,15){\line(1,0){7}}
\put(29,20){\line(0,1){13}}
\qbezier(29,20)(29,18.5)(30.5,18.5)
\put(29,36){\line(0,1){5}}

\put(30,31){\line(1,0){8}}
\qbezier(30.5,18.5)(32,18.5)(32,17)
\qbezier(34,35)(37,35)(37,32)
\put(34,46){\line(1,0){8}}
\put(34,50){\line(1,0){8}}
\put(37,10){\line(0,1){20}}
\qbezier(37,10)(37,8)(39,8)
\put(39,11){\line(1,0){24}}
\put(39,15){\line(1,0){24}}
\qbezier(38,31)(41,31)(41,28)
\qbezier(39,8)(41,8)(41,10)

\put(41,16){\line(0,1){12}}
\qbezier(42,46)(47,46)(47,41)
\qbezier(42,50)(51,50)(51,41)
\qbezier(47,41)(47,33)(55,33)

\qbezier(51,41)(51,37)(55,37)
\put(55,33){\line(1,0){8}}
\put(55,37){\line(1,0){8}}
\qbezier(58,32)(58,24)(66,24)
\put(58,38){\line(0,1){4}}
\qbezier(58,42)(58,50)(66,50)

\qbezier(62,32)(62,28)(66,28)
\put(62,38){\line(0,1){4}}
\qbezier(62,42)(62,46)(66,46)
\qbezier(63,11)(71,11)(71,19)
\qbezier(63,15)(67,15)(67,19)
\qbezier(63,33)(67,33)(67,29)
\qbezier(63,37)(71,37)(71,29)
\put(66,24){\line(1,0){8}}
\put(66,28){\line(1,0){8}}
\put(66,46){\line(1,0){24}}
\put(66,50){\line(1,0){24}}
\put(67,19){\line(0,1){4}}

\put(71,19){\line(0,1){4}}
\qbezier(74,24)(78,24)(78,20)
\qbezier(74,28)(82,28)(82,20)
\qbezier(78,20)(78,11)(87,11)

\qbezier(82,20)(82,15))(87,15)
\put(87,11){\line(1,0){8}}
\put(87,15){\line(1,0){8}}
\put(88,33){\line(0,1){12}}
\qbezier(88,33)(88,30)(91,30)
\qbezier(88,51)(88,53)(90,53)

\qbezier(90,53)(92,53)(92,51)
\put(91,30){\line(1,0){8}}
\qbezier(92,29)(92,26)(95,26)
\put(92,31){\line(0,1){20}}
\put(93,46){\line(1,0){8}}
\put(93,50){\line(1,0){8}}
\qbezier(95,11)(104,11)(104,20)
\qbezier(95,15)(100,15)(100,20)
\put(95,26){\line(1,0){10}}
\qbezier(97,44)(97,42.5)(98.5,42.5)
\qbezier(97,51)(97,54.5)(100.5,54.5)
\qbezier(98.5,42.5)(100,42.5)(100,41)

\put(100,20){\line(0,1){5}}
\put(100,28){\line(0,1){13}}
\qbezier(100.5,54.5)(104,54.5)(104,51)
\put(104,20){\line(0,1){5}}
\put(104,28){\line(0,1){23}}
\qbezier(105,26)(114,26)(114,35)
\qbezier(105,30)(110,30)(110,35)
\qbezier(105,46)(110,46)(110,41)
\qbezier(105,50)(114,50)(114,41)
\put(110,35){\line(0,1){6}}
\put(114,35){\line(0,1){6}}

\put(30, -1){Figure 1. A smoothly freely slice link.}
\end{picture}

The ribbon link $L$ of Figure 1 is a satellite of the Hopf link,
formed  by replacing each component by a copy of $KT$.
Let $Z_L$ be the exterior of the pair of slice discs obtained from the obvious ribbon discs.
Then $\pi=\pi_1Z_L\cong{F(2)}$ (see  \cite[Chapter 1.7]{HiA}) and $Z_L$ is aspherical,
as in Theorem \ref{Wh(L)}.
The torus which separates the two components of the Hopf link splits 
$M(L)=\partial{Z_L})$ into two homeomorphic parts $U$ and $V$, 
each of which is the exterior of a 2-component link,
with one component the knot $KT$ and the other a meridian $\mu_{KT}$ of the knot.
If $\pi_2(U)\not=0$ then there is an essential 2-sphere $\Sigma\subset{U}$,
by the Sphere Theorem.
This must bound a  ball $B$ in $M(KT)$, since $M(KT)$ is aspherical.
If $\sigma$ is essential in $U$ then $B\not\subset{U}$, and so $\mu_{KT}\subset{B}$.
But this contradicts the fact the image of $\mu_{KT}$ generates $H_1(M(K);\Z)\cong\Z$.
Therefore $U$ (and hence $V$) is aspherical, and the inclusion of $T\subset{U}$ is $\pi_1$-injective.
Hence $M(L)$ is aspherical also.

\section{cohomological dimension $2$}

If $\cd\pi=2$ and $w:\pi\to\Z^\times$ is a homomorphism then there is a $PD_4$-pair $(X,Y)$ 
with $X\simeq{K(\pi,1)}$ if and only if there is a map $p:Y\to{K(\pi,1)}$ 
such that $\pi_1p$ is an epimorphism,  $H^q(p;\Z\pi)$ is an isomorphism
for $q\leq2$ and $w_1Y=wp$, by Theorems \ref{cd<n} and \ref{PD4pair}.
We shall show that the study of such $PD_4$-pairs may largely 
be reduced to the cases with $Y$ aspherical.
In the latter part of this section we shall focus on the 4-manifold case.
One significant feature is that if $\cd\pi=2$ then $H^2(\pi;\Z\pi)\not=0$,
and so no aspherical $PD_4$-pair realizing $\pi$ has $\pi_1$-injective boundary.

If $M$ is a compact aspherical smooth 4-manifold such that $\pi_1\partial{M}$ maps onto $\pi=\pi_1M$ then 
we may use handle trading to show that $M$ may be obtained from $\partial{M}\times[0,1]$
by adding 3- and 4-handles only.
Inverting the handle decomposition, 
we see that $M$ is homotopy equivalent to a finite 2-complex.
Hence if also $\cd\pi=2$ then $\gd\pi=2$.
Conversely, 
every finite 2-complex is homotopy equivalent to a compact smooth 4-manifold with boundary,
and so every group $\pi$ with $\gd\pi\leq2$ is the fundamental group of a compact aspherical 4-manifold.

If $G$ is a group and $\Gamma$ is a right $\Z{G}$-module let
\[
H_q(G,1;\Gamma)=H_q(K(G,1),*;\Gamma),
\]
for all $q$. 
There is then a 4-term exact sequence of homology 
\[
0\to{H_1(G;\Gamma)}\to{H_1(G,1;\Gamma)}\to{H_0(*;\Gamma)}\to{H_0(G;\Gamma)}\to0.
\]
(When $G$ is a knot group and $\Gamma=\Z{G^{ab}}\cong\Z[t,t^{-1}]$ 
this sequence is well known to knot theorists,
for then $H_1(G;\Z{G^{ab}})=G'/G''$ and $H_1(G,1;\Z[t,t^{-1}])$ is the Alexander module
for the knot group.)

\begin{lemma}
\label{Kabproj}
Let $f:A*B\to\pi$ be  an epimorphism and let $K=\ker(f)$.
If $\cd \pi\leq2$ then $K^{ab}$ is a projective $\Z\pi$-module.
\end{lemma}

\begin{proof}
If we set  $\Gamma=\Z{G}$ in the above 4-term exact sequence
we see that $H_1(G,1;\Z{G})\cong\mathcal{I}(G)$,
the augmentation ideal in $\Z{G}$.
Taking $G=A*B$ and $\Gamma=\Z\pi$ (with the right $\Z{G}$-module structure determined by $f$),
we get the sequence
\[
0\to{K^{ab}}\to(\Z\pi\otimes_A\mathcal{I}(A))\oplus(\Z\pi\otimes_B\mathcal{I}(B))\to\Z\pi\to\Z\to0.
\]
If $\cd\pi\leq2$ then  $\cd{A}\leq2$ and  $\cd{B}\leq2$ also,
and so $\mathcal{I}(A)$ and $\mathcal{I}(B)$ have projective dimension $\leq1$ over 
$\Z{A}$ and $\Z{B}$,  respectively. 
Thus there is an exact sequence
\[
0\to{P_2}\to{P_1}\to(\Z\pi\otimes_A\mathcal{I}(A))\oplus(\Z\pi\otimes_B\mathcal{I}(B))\to0
\]
with $P_1$ and $P_2$ projective $\Z\pi$-modules, and so there is an exact sequence
\[
0\to{K^{ab}}\oplus{P_2}\to{P_1}\to\Z\pi\to\Z\to0.
\]
The left hand module ${K^{ab}}\oplus{P_2}$ is projective, by Schanuel's lemma, since $\cd\pi\leq2$.
Hence $K^{ab}$ is projective.
\end{proof}

\begin{lemma}
\label{ddual}
Let $G$ be an $FP$ group with $\cd{G}=n$. 
\begin{enumerate}
\item{}The right $\Z{G}$-module $H^n(G;\Z{G})$ has no non-trivial projective quotient.
\item{}If $G$ is a duality group then $H^n(G;\Z{G})$ has no proper direct summand.
\end{enumerate} 
\end{lemma}

\begin{proof}
Since $G$ is $FP$ the augmentation module $\Z$ has a  resolution $P_*\to\Z$ 
by finitely generated projective left $\Z{G}$-modules,
with $P_i=0$ for $i>n$.
Let $P^i=P_i^*=\Hom_{\Z{G}}(P_i,\Z{G})$ be the dual complex of projective right $\Z{G}$-modules.
Then $H^n(G;\Z{G})$ is a quotient of $P^n$.
If $H^n(G;\Z{G})$ has a projective quotient $Q$ then the composite projection 
$P^n\to{H^n(G;\Z{G})}\to{Q}$ splits,
and so $P^n\cong{P\oplus{Q}}$, with $\mathrm{Im}(\delta^n)\leq{P}$.
Hence $P_n=P^{n*}\cong{P^*}\oplus{Q^*}$.
Let $\delta_P$ and $\delta_Q$ be the restrictions of $\delta_n$ to these summands.
Then $\delta_n=\delta_P\oplus\delta_Q$, 
and $\delta^n:P^{n-1}\to{P^n}$ is the transpose of $\delta_n$.
Since $\delta^n$ has image in $P$, we see that $\delta_Q=0$.
Hence $Q=0$, since  the $n$th differential $\delta_n$ is a monomorphism. 

Suppose now that $G$ is a duality group.
The dual complex $P^*$ is then a finite projective resolution for $H^n(G;\Z{G})$,
and dualizing again recovers the resolution $P_*$.
If $H^n(G;\Z{G})\cong{A\oplus{B}}$, where $A\not=0$
then $\Ext^n_{\Z{G}}(A\oplus{B},\Z{G})\cong\Z$,
while $\Ext^q_{\Z{G}}(A\oplus{B},\Z{G})=0$ for all $q\not=n$.
Since $\Z$ is not a proper direct sum, this is only possible if
$\Ext^q_{\Z{G}}(B,\Z{G})=0$ for all $q$, in which case $B=0$.
\end{proof}

Our main result in this section is the following theorem,
which is essentially Theorem \ref{B} of the introduction.

\begin{theorem}
\label{cd2asph}
Let $(X,\partial{X})$ be an aspherical $PD_4$-pair and let $\pi=\pi_1X$.
If $\cd\pi=2$ then there is a set of $PD_4$-pairs  $\{(X_1,N_1),\dots,(X_k,N_k)\}$, 
where $X_i$ is aspherical and $N_i$ is aspherical or is a copy of $S^2\times{S^1}$ 
or $S^2\tilde\times{S^1}$,
for $i\leq{k}$, 
and a $PD_3$-homology sphere $\Sigma$ such that
$(X,\partial{X})\simeq(\natural_{i=1}^k(X_i,N_i))\natural(C,\Sigma)$, 
where $C$ is the cone on $\Sigma$.
\end{theorem} 

\begin{proof}
The boundary $\partial{X}$ is connected and $\pi_1\mathrm{inc}$ is an epimorphism,
since $\cd\pi=2$,
and so $\partial{X}$ has an indecomposable summand $N_1$ 
such that $\nu_1=\pi_1N_1$ maps nontrivially to $\pi$. 
Let $\omega$ be the free product of all of the other prime factors of $\pi_1\partial{X}$. 
Then $\pi_1\partial{X}\cong\nu_1*\omega$.
Let $\rho_1$ and $B$ be the images of $\nu_1$ and $\omega$ in $\pi$.
Then $\rho_1\not=1$ and the kernel of the map from $\rho_1*B$ to $\pi$ induced by $\pi_1\inc$
is a free group $K$ such that $K^{ab}$ is a projective $\Z\pi$-module, by Lemma \ref{Kabproj}.
This kernel is also a quotient of $\kappa=\ker(\pi_1\inc)$,
which has abelianization isomorphic to 
$H_1(\partial{X};\Z\pi)\cong\overline{H^2(\pi;\Z\pi)}=\overline{\Ext^2_{\Z\pi}(\Z,\Z\pi)}$.
The $\Z\pi$-module $H^2(\pi;\Z\pi)$ has no projective quotient,  by Lemma \ref{ddual},
so we must have $K^{ab}=1$.
Since $K$ is free, it follows that $K=1$, and so $\pi\cong{\rho_1*B}$.

A finite induction then shows that $\partial{X}\cong{N\#\Sigma}$,
where $N$ is a connected sum of $k$ indecomposable $PD_3$-complexes $N_i$
with $\nu_i=\pi_1N_i$ an extension of a factor $\rho_i$ of $\pi$ by a normal subgroup,
$\pi=*_{i=1}^k\rho_i$ and $\sigma=\pi_1\Sigma$ maps trivially to $\pi$.
Each summand $N_i$ is aspherical, or a copy of $S^2\times{S^1}$ or $S^2\tilde\times{S^1}$, 
by Theorem \ref{asph4mfdbdry}.
Now $H_1(\partial{X};\Z\pi)\cong{H^2(\pi;\Z\pi)}$ splits as the direct sum 
\[
\oplus_{i=1}^kH_1(\nu_i;\Z\pi)\oplus{H_1(\sigma;\Z\pi)}
\cong
\oplus_{i=1}^kH_1(\nu_i;\Z\pi)\oplus(\sigma^{ab}\otimes_\Z\Z\pi),
\]
since $\sigma\leq\kappa$.
The module $H^2(\pi;\Z\pi)$ is torsion free as an abelian group
and has no free direct summand as a right $\Z\pi$-module.
Hence $\sigma^{ab}=0$ and so $\Sigma$ is a homology 3-sphere.

Let $X_i=MCyl(c_i)$ be the mapping cylinder of a map $c_i:N_i\to{K(\rho_i,1)}$ 
corresponding to the epimorphism $p_i:\nu_i\to\rho_i$ induced by $\inc$.
Since $H^q(p;\Z\pi)=\oplus_{i=1}^kH^q(p_i;\Z\pi)$ for all $q$,
each term $H^q(p_i;\Z\pi)$ is an isomorphism for $q\leq2$.
Hence $H^q(p_i;\Z\rho_i)$ is an isomorphism for $q\leq2$,
and so $(X_i,N_i)$ is a $PD_4$-pair, 
by Theorem \ref{PD4pair}.
Let $(X',\partial{X'})=\natural_{i=1}^k(X_i,N_i)$ and let
$C=c\Sigma$ be the cone on $\Sigma$.
Then $\partial{X'}\simeq{N}$, and there are isomorphisms $s:\pi\cong*\rho_i$ 
and $b:\pi_1N\cong*\nu_i$ such that $sf=(*c_i)b$.
These may be realized by a map from $(X',\partial{X'})\natural(C,\Sigma)$
to $(X,\partial{X})$ which is a homotopy equivalence of pairs.
\end{proof}

\begin{corollary}
\label{gd2oneend}
Let $(M,\partial M)$ be a compact aspherical $4$-manifold with $\pi = \pi_1M$.
If $\gd\pi=2$ and $\pi$ has one end then there is a $4$-manifold $M'\simeq{M}$
such that $N=\partial{M'}$ is aspherical and a contractible $4$-manifold $C$ such that 
$(M,\partial{M})\simeq(M',N)\natural(C,\partial{C})$.
\qed
\end{corollary}

The summand $N$ of $\partial{M}$ may itself be a homology 3-sphere. 
For instance, 
if $\pi$ is the Higman group, with presentation
\[
\langle{a,b,c,d}\mid{bab^{-1}=a^2},~cbc^{-1}=b^2,~dcd^{-1}=c^2,~ada^{-1}=a^2\rangle,
\]
then the associated 2-complex is aspherical,  so $\gd\pi=2$,  and $H_p(\pi;\Z)=0$ for $p>0$.
In this case $\partial{M}$ is a homology 3-sphere.

For the remainder of this section we shall assume that $\pi$ is a duality group of dimension 2.
Let  $\mathcal{D}_\pi=H^2(\pi;\Z\pi)=\Ext^2_{\Z\pi}(\Z,\Z\pi)$ be the dualizing module.
(This is naturally a right $\Z\pi$-module.)
If $\pi$ is a $PD_2$-group then $\mathcal{D}_\pi\cong\Z$
(as an abelian group) and every extension of $\pi$ by 
$\overline{\mathcal{D}_\pi}$ gives a $PD_3$-group.

\begin{theorem}
\label{PD2cor}
Let $(X,\partial{X})$ be an aspherical $PD_4$-pair such that $\pi = \pi_1X\cong\pi_1F$,
where $F$ is an aspherical closed surface, and suppose that the boundary components are
closed $3$-manifolds.
Then $(X,\partial{X})\simeq(E\cup_{\partial{E}}MCyl(f),N)\natural(C,\partial{C})$,
where $E$ is the total space of a $D^2$-bundle over $F$,
$N$ is a compact aspherical $3$-manifold, 
$f:N\to\partial{E}$ is a degree-$1$ map which is a $\Z\pi$-homology equivalence 
and $C$ is a contractible $4$-manifold.
\end{theorem}

\begin{proof}
We may assume that $(M,\partial{M})\simeq(X',N)\natural(C,\partial{C})$,
where $X'\simeq{F}$, $N$ is a prime 3-manifold and $C$ is contractible,
by Theorem \ref{cd2asph}.

Let $j_*:\nu=\pi_1N\to\pi$ be the homomorphism induced by the inclusion $j:N\to{X'}$.
Then $j_*$ is an epimorphism, and $\kappa=\mathrm{Ker}(j_*)$ has abelianization
$\kappa^{ab}\cong\overline{H^2(\pi;\Z\pi)}$. 
Thus $\kappa^{ab}\cong\Z$ as an abelian group,
and $\nu/\kappa'$ is an extension of $\pi$ by $\Z$.
In particular $N$ is aspherical, since $\nu\not\cong\Z$.
Moreover, $H_2(N;\Z\pi)=0$, 
since $H_2(X';\Z\pi)=H_2(\pi;\Z\pi)=0$ and $H_3(X',N;\Z\pi)\cong{H^1(\pi;\Z\pi)}=0$.

The group $\nu/\kappa'$ is the fundamental group of the total space
of an $S^1$-bundle over $F$.
Let $E$ be the total space of the associated disc bundle.
The projection of $\nu$ onto $\nu/\kappa'$ is realized by a map 
$f:N\to\partial{E}$, since $\partial{E}$ is aspherical.
Now $H_2(N;\Z\pi)=H^1(N;\Z\pi)=0$, by part (1) of Theorem \ref{cd<n}.
Comparison of the spectral sequences for the projections of $N$ and $\partial{E}$
onto $F$ shows that $f$ is a $\Z\pi$-homology equivalence.
Hence $(X',N)\simeq(E\cup_{\partial{E}}MCyl(f),N)$, 
by a Mayer-Vietoris argument, and so
$(X,\partial{X})\simeq(E\cup_{\partial{E}}MCyl(f),N)\natural(C,\partial{C})$.
\end{proof}

The map $f$ need not be a $\Z\pi_1\partial{E}$-homology equivalence.
For instance, suppose that $E=T_g\times{D^2}$ and replace a neighbourhood $D^2\times{S^1}$ of a fibre
$S^1$ in $\partial{E}=T_g\times{S^1}$ by the exterior  $X(K)$ of a nontrivial knot $K$,
aligning the meridian and longitude in $\partial{X(K)}$ with the fibre $S^1$ and $\partial{D^2}$, respectively.
The resulting closed 3-manifold is aspherical and maps to $\partial{E}$ via a $\Z\pi_1T_g$-homology equivalence. If $K$ has non-trivial Alexander polynomial this is not
a $\Z\pi_1\partial{E}$-homology equivalence.

There are no obvious ``basic" examples if $\pi$ is not a $PD_2$-group.
In view of \cite{dh2}, the case of greatest interest is when $\pi$ is  solvable,
and thus is a Baumslag-Solitar group $BS(1,m)$, for some $m\not=0$.

It can be shown that if $\pi$ is not a $PD_2$-group then $\kappa$ must have trivial centre.
(For the centre $\zeta\kappa$ is normal in $\pi_1\partial{X}$, 
and so if $\zeta\kappa\not=1$, 
then $\partial{X}$ would be homotopy equivalent to a Seifert fibred 3-manifold. 
This leads to a contradiction with the fact that $H_2(\kappa;\Z)=0$.)

If $\kappa$ has trivial centre then extensions of $\pi$ by $\kappa$ correspond
to homomorphisms $\alpha:\pi\to\mathrm{Out}(\kappa)$ \cite[Corollary 4.6.8]{Brown},
and there are infinitely many extensions realizing 
$\alpha^{ab}:\pi\to\mathrm{Aut}(\kappa^{ab})$,
but even if one such is a $PD_3$-group,
it is not clear that other extensions inducing the same automorphism of $\kappa^{ab}$ must
give $PD_3$-groups.

The example following Theorem \ref{PD2cor} shows that 
$\kappa$ need not be a free group, 
even if $\partial{X}$ is aspherical.
However, 
if $\kappa$ is a free group then $(X,\partial{X})$ is minimal with respect to the partial order 
determined by degree-1 maps $f:(X',\partial{X'})\to(X,\partial{X})$ such that $\pi_1f$ is an isomorphism.
This follows from the facts that $H_2(\partial{X};\Z\pi)={H^1(\partial{X};\Z\pi)=0}$, 
by part (1) of Theorem \ref{cd<n}, 
the dualizing module $\mathcal{D}_\pi$ is Hopfian as a $\Z\pi$-module
(since $End_{\Z\pi}(\mathcal{D}_\pi)\cong{End_{\Z\pi}(\Z)}=\Z$) and the main result of 
\cite{St65}.
If there is such an aspherical $PD_4$-pair $(X,\partial{X})$ with $\kappa$ a free group
then $D_\pi$ is a free abelian group, since $\kappa^{ab}\cong\overline{\mathcal{D}_\pi}$.
Whether this is so for all duality groups of dimension 2 is an open question.

It is not clear whether for every pair $(\pi,w)$, 
with $\pi$ a duality group of dimension 2 and $w:\pi\to\Z^\times$ a homomorphism
there is a $PD_3$-group $G$ and an epimorphism $\theta:G\to\pi$ such that
$w_1G=w\theta$,  $\mathrm{Ker}(\theta)$ is free and condition (2) of Theorem \ref{PD4pair} holds.
If $E=\Hom_\Z(\overline{\mathcal{D}_\pi},\Z\pi)$ has the left $\Z\pi$-module structure 
determined by $(g.\lambda)(d)=g\lambda(g^{-1}d)$,
for all $\lambda\in{E}$ and $d\in\overline{\mathcal{D}_\pi}$, 
then the LHS homology spectral sequence with coefficients $\Z\pi$ for $G$ as an extension of $\pi$ 
and Poincar\'e duality together imply that $H^1(\pi;E)=0$ and $H^2(\pi;E)\cong\Z$.
(Note that $H^0(\pi;E)=\Hom_{\Z\pi}(\overline{\mathcal{D}_\pi},\Z\pi)=0$.)
However it is not easy to test these conditions,
and we do not know whether they always hold,
nor whether they imply that there is such a group $G$.

As in the case of free groups, 
there are examples with $\pi$ a proper free product and $\partial{M}$ aspherical,
and so we cannot always reduce to the case when $\pi$ is  indecomposable.
Perhaps the simplest such is the one with $\pi=\pi_1M\cong\Z^2*\Z^2$
and the Kirby calculus presentation given in Figure 2,
which is based on a similar presentation for $(T^2\times{D^2})\natural(T^2\times{D^2})$.
The dotted components represent 1-handles (and generators for $\pi$),
while the other two components represent attaching maps of 2-handles
(and relators for $\pi$).

\setlength{\unitlength}{1mm}
\begin{picture}(90,55)(4,-1)

\put(12,14.1){$\bullet$}
\put(9,11){$a$}

\put(5,20){\line(0,1){13}}
\qbezier(5,33)(5,38)(11,38)
\qbezier(5,20)(5,15)(10,15)

\put(10,15){\line(1,0){5}}
\put(10,38){\line(1,0){2}}
\put(13,27){\line(0,1){14}}
\qbezier(13,27)(13,21)(19,21)
\qbezier(13,41)(13,46)(18,46)
\put(14,38){\line(1,0){1}}
\qbezier(15,15)(20,15)(20,20)
\qbezier(15,38)(20,38)(20,33)
\put(18,46){\line(1,0){40}}
\put(19,29){\line(1,0){17}}
\qbezier(19,29)(16.5,29)(16.5,26.5)
\qbezier(16.5,26.5)(16.5,24)(19,24)

\put(20,20){\line(0,1){8}}
\put(20,30){\line(0,1){3}}
\put(21,21){\line(1,0){17}}
\put(21,24){\line(1,0){15}}
\put(28,18){$e$}

\put(37,22){\line(0,1){11}}
\qbezier(37,33)(37,38)(42,38)
\qbezier(37,20)(37,15)(42,15)
\qbezier(38,21)(39.5,21)(39.5,22.5)
\qbezier(38,24)(39.5,24)(39.5,22.5)
\put(38,29){\line(1,0){22}}

\put(41,11){$k$}
\put(42,15){\line(1,0){5}}
\put(42,38){\line(1,0){5}}
\put(44,14.1){$\bullet$}
\qbezier(47,15)(52,15)(52,20)
\qbezier(47,38)(52,38)(52,33)

\put(52,20){\line(0,1){8}}
\put(52,30){\line(0,1){3}}
\qbezier(58,35)(58,32)(61,32)
\put(58,35){\line(0,1){4}}
\qbezier(58,39)(58,42)(61,42)
\qbezier(58,46)(62,46)(62,42)

\qbezier(60,29)(62,29)(62,31)
\put(61,32){\line(1,0){28}}
\put(62,34){\line(0,1){8}}
\put(63,42){\line(1,0){49}}

\put(75,20){\line(0,1){11}}
\qbezier(75,33)(75,38)(80,38)
\qbezier(75,20)(75,15)(80,15)
\put(79,11){$w$}

\put(80,38){\line(1,0){5}}
\put(80,15){\line(1,0){5}}
\put(82,14.1){$\bullet$}
\qbezier(85,38)(90,38)(90,33)
\qbezier(85,15)(90,15)(90,20)
\qbezier(87.5,22.5)(87.5,21)(89,21)
\qbezier(87.5,22.5)(87.5,24)(89,24)
\put(89,21){\line(1,0){17}}

\put(90,22){\line(0,1){11}}
\put(91,24){\line(1,0){15}}
\put(91,32){\line(1,0){17}}
\put(97.5,18){$s$}

\put(107,20){\line(0,1){11}}
\qbezier(107,33)(107,38)(112,38)
\qbezier(107,20)(107,15)(112,15)
\qbezier(108,32)(113,32)(113,28)
\qbezier(108,24)(113,24)(113,28)
\put(108,21){\line(1,0){3}}

\qbezier(111,21)(116,21)(116,26)
\qbezier(111,42)(116,42)(116,37)
\put(112,11){$z$}
\put(112,15){\line(1,0){5}}
\put(112,38){\line(1,0){3}}
\put(115,14.1){$\bullet$}
\put(116,26){\line(0,1){11}}
\qbezier(117,38)(122,38)(122,33)
\qbezier(117,15)(122,15)(122,20)
\put(122,20){\line(0,1){13}}

\put(30, 3){Figure 2.  $\pi=\Z^2*\Z^2$ and $\partial{M}$ aspherical.}

\end{picture}

It is easy to see that $M\simeq{T^2}\vee{T^2}$,  and so $M$ is aspherical.
Adding relators representing the longitudes of each component to a presentation of the link group
gives the presentation
\[
\langle{a,e,k,s,w,z}\mid {ae=ea},~ek=ke,~sw=ws,~sz=zs,
\]
\[
e=[z,w],~s=[k,a]\rangle
\]
for $\nu=\pi_1\partial{M}$.
Since the relations imply that $es=se$, 
this group is a generalized free product of two copies of $F(2)\times\Z$ 
with amalgamation over $\Z^2$, and so it has one end. 
Hence $\partial{M}$ is aspherical.
(In fact $\nu=\pi_1N$, where $N$
 is the union of two copies of $T_o\times{S^1}$ along a torus.)

\section{3-dimensional duality groups}

When $\cd\pi=3$ the boundary need not be connected, 
and the images of the fundamental groups of boundary components 
may have infinite index in $\pi$.
Corollary \ref{3DGcor} below is an analogue of Theorem \ref{cd2asph} 
for the case when $\pi$ is a 3-dimensional duality group.

\begin{lemma}
\label{PD4pairs of groups}
Let $(X,\partial{X})$ be an aspherical $PD_4$-pair and let $\pi=\pi_1X$.
If the boundary components  are $\pi_1$-injective 
then the following are equivalent;
\begin{enumerate}
\item{} the boundary components are aspherical;
\item$\pi$ has one end;
\item$\pi$ is a duality group of dimension $3$.
\end{enumerate}
If these conditions hold then $(X,\partial{X})$ is a $PD_4$-pair of groups.
\end{lemma}

\begin{proof}
The equivalence of these conditions follows from consideration of the homology
of the pair $(K(\pi,1),\partial{X})$ with coefficients $\Z\pi$, together with the fact that
since $H^3(\pi;\Z\pi)\cong{H_1(X,\partial{X};\Z\pi)}=\widetilde{H}_0(\partial{X};\Z\pi)$,
it is free as an abelian group.
\end{proof}

If $\pi$ is a duality group  of dimension 3 but is not a $PD_3$-group 
then $H^3(\pi;\Z\pi)$ is not finitely generated \cite{Fa75},
and so for each boundary component $\partial_jX$
the image of $\pi_1\partial_jX$ in $\pi$ has infinite index.

Let $\mathcal{G}=\{\gamma_i:i\leq{k}\}$ be a finite set of monomorphisms  $\gamma_i:G_i\to\pi$,
and let $(K(\pi,\mathcal{G}),\sqcup{K_i})$ be the pair with ambient space the union of $K(\pi,1)$ with the mapping cylinders of the maps from $K_i=K(G_i,1)$ to $K(\pi,1)$ determined by the $\gamma_i$s.
If this pair is a $PD_4$-pair we shall write also just $(\pi,\mathcal{G})$.
(In our considerations of $PD_4$-pairs of groups 
we want to allow pairs such as $(K\times[0,1],K\times\{0,1\})$ with $K=K(\pi,1)$,
and so we label the inclusions,  rather than the images of the groups $G_i$.)

The next result includes part (1) of Theorem \ref{C} of the introduction.

\begin{theorem} 
\label{PD4duality}
Let $(X,\partial{X})$ be an aspherical $PD_4$-pair such that
$\pi=\pi_1X$ is  a duality group of dimension $3$ and $\partial{X}=\sqcup_{i=1}^mY_i$.
Let $G_i$ be the image of $\pi_1Y_i$ in $\pi$, for $i\leq{m}$, 
and let $\mathcal{K}$ be the set of inclusions of the $G_i$ in $\pi$.
Then
\begin{enumerate}
\item$Y_i=Y_i'\#\Sigma_i$, where $Y_i'$ is aspherical,
$\Sigma_i$ is a $PD_3$-homology sphere
and $G_i$ is a $PD_3$-group, for all $i\leq{m}$;
\item$(X,\partial{X})\simeq(X',Y')\natural(C_1,\Sigma_1)\dots\natural(C_m,\Sigma_m)$, 
where $(X',\partial{X'})$ is a $PD_4$-pair, 
$X'\simeq{X}$, $\partial{X'}=\sqcup_{i=1}^mY_i'$, 
$C_i$ is  the cone over $\Sigma_i$
and $C_i\cap\partial{X'}=\overline{\Sigma_i\setminus{D^3}}\subset{Y_i'}$, for all $i\leq{m}$;
\item$(\pi,\mathcal{K})$ is a $PD_4$-pair of groups and
there are degree-$1$ maps $p_j:Y_i'\to{K_i}$ with $\mathrm{Ker}(\pi_1p_i)$ acyclic,
for $i\leq{m}$,
such that\\ 
 $(X',\partial{X'})\simeq(\pi,\mathcal{K})\cup_{K_1}(MCyl(p_1),Y_1')\dots
\cup_{K_m}(MCyl(p_m),Y_m')$.
\end{enumerate}
\end{theorem} 

\begin{proof}
The subgroups $G_i$ are $PD_3$-groups, by part (6) of Theorem \ref{cd<n}
and so are indecomposable as free products.
Since $\ker\pi_1\inc_i$ is perfect,
 $Y_i\cong{Y_i'}\#\Sigma_i$,  where $Y_i'$ is indecomposable,  
$\pi_1Y_i'$ maps onto $G_i$ and $\Sigma_i$ is a $PD_3$ homology 3-sphere 
such that $\pi_1\Sigma_i\leq\ker\pi_1\inc_i$,
by Lemma \ref{G*Hontopi}.
Since $\pi_1Y_i'$ maps onto a $PD_3$-group which is a subgroup of $G_i$,
it must be aspherical, by Theorem \ref{asph4mfdbdry}.
(Note that indecomposable $PD_3$-homology spheres are either aspherical or have
finite fundamental group $I^*=SL(2,5)$ \cite{Hi20}.)

Let $f_i:Y_i'\to{K_i}$ be the map which induces the composite of the inclusion of 
$\pi_1Y_i'$ into $\pi_1N_i$ with $\pi_1\mathrm{inc}_i$.
A Mayer-Vietoris argument applied to $K(\pi_1Y_i,1)\simeq{Y_i'}\vee{K(\pi_1\Sigma_i,1)}$,
shows that $\mathrm{Ker}(\pi_1f_i)$ is acyclic, and it follows easily that $f_i$ has degree 1.
It now follows from Lemma \ref{excismv} that $(\pi,\mathcal{K})$ is a $PD_4$-pair of groups.
The rest of the argument is as before.
\end{proof}

Thus $PD_4$-groups with boundary and ambient group a 3-dimensional duality group 
may be obtained from $PD_4$-pairs of groups by adding mapping cylinders of degree-1 maps 
which satisfy the above acyclicity condition, 
together with boundary connected sums with contractible $PD_4$-complexes.

\begin{corollary}
\label{3DGcor}
Let $\pi$ be a duality group of dimension $3$. Then
\begin{enumerate}
\item{}$\pi=\pi_1X$ for some  aspherical $PD_4$-pair $(X,\partial{X})$ 
if and only if $\pi$ is the ambient group of a $PD_4$-pair of groups $(\pi,\mathcal{K})$;
\item{}Any such $PD_4$-pair $(X,\partial{X})$ with $m$ boundary components $Y_i$
is homotopy equivalent to
\[
(\pi,\mathcal{K})\cup_{K_1}(M_1,K_1\sqcup{Y_1})\natural(C_1,\Sigma_1)\dots
\cup_{K_m}(M_m,K_m\sqcup{Y_m})\natural(C_m,\Sigma_m),
\]
where $\mathcal{K}$ is determined by the inclusion of $Y$ into $X$,
$M_i$ is the mapping cylinder of a $\Z{G_i}$-homology equivalence from $Y_i$ to $K_i$
and $C_i$ is the cone over a $PD_3$-homology sphere $\Sigma_i$.
\qed
\end{enumerate}
\end{corollary}

It seems unlikely that the duality groups of dimension 3 of interest here
should have a simpler characterization than that in part (1) of the Corollary.
Let $T_o$ be the punctured torus and $D_{oo}$ be the twice punctured disc.
Then $T_o\times{T}$ and $D_{oo}\times{T}$ are 4-manifolds with fundamental group $F(2)\times\Z^2$.
Thus the ambient group does not determine the boundary, 
and need not be a $PD_3$-group.  
On the other hand,  the groups $BS(1,m)\times\Z$ 
are solvable duality groups of dimension 3 with finite $3$-dimensional $K(\pi,1)$ complexes.
If $m\not=\pm1$ then $BS(1,m)\times\Z$ has no subgroups which are $PD_3$-groups,
and so is not realizable by an aspherical $PD_4$-pair $(X,\partial{X})$.

If $\pi=\pi_1N$ for some closed 3-manifold $N$ then clearly $\pi\cong\pi_1M$, 
where $M=N\times[0,1]$.
However it is not known whether every $PD_3$-group is the fundamental group 
of a compact aspherical 4-manifold,
let alone whether every $PD_3$-pair of groups is realized by a compact aspherical 3-manifold with boundary. 
Thus it would be unreasonable to expect to prove now that 
$(K(\pi,\mathcal{K}),\sqcup{K_i})$ is homotopy equivalent to a compact 4-manifold.
Beyond this, there is the general problem of surgery,
which at the time of writing cannot be used to prove
that a compact aspherical 4-manifold is homeomorphic to one assembled 
as in Corollary \ref{3DGcor}.

\section{the other groups with $H^2(\pi;\Z\pi)=0$}

We show next that every aspherical $PD_4$-pair with $H^2(\pi;\Z\pi)=0$ 
may be obtained from pairs with $\pi_1$-injective boundary
by adding mapping cylinders over boundary components,
as in Theorem \ref{PD4duality} and Corollary \ref{3DGcor}.
The case when $\pi$ is a free group is covered by Corollary \ref{freeboundaryinj},
and so we may assume that $\cd\pi=3$.

\begin{lemma}
\label{end module}
Let $\pi=(*_{i=1}^mG_i)*F(r)$, where each factor $G_i$ is $FP_2$ and has one end.
If $\pi\not=1$ then
\begin{enumerate}
\item{}if $\pi\cong{F(r)}$ then $H^1(\pi;\Z\pi)$ has projective dimension $1$;
\item{}if $m>0$, so that $\pi$ is not free, then $H^1(\pi;\Z\pi)\cong(\Z\pi)^{m+r-1}$.
 \end{enumerate}
\end{lemma}

\begin{proof}
Let $(\mathcal{G},\Gamma)$ be a finite graph of groups,
with all edge groups finite and all vertex groups finite or one-ended, and let $\pi=\pi\mathcal{G}$.
Let $E$ be the set of edges and $V$  the set of vertices in the graph $\Gamma$.
If $\pi$ is infinite then $H^0(\pi;\Z\pi)=0$,
and a Mayer-Vietoris argument due to Chiswell \cite{Ch76} gives
an exact sequence
\[
0\to\oplus_{v\in{V}}H^0(G_v;\Z\pi)\to
\oplus_{e\in{E}}H^0(G_e;\Z\pi)\to{H^1(\pi;\Z[\pi])}\to0.
\]
The sequence is exact at the right since the vertex groups are finite or have one end, 
so the next term $\oplus_{v\in{V}}H^1(G_v;\Z\pi)$ is $0$.

If $\pi=F(r)$ then we may assume that $\Gamma$ has one vertex and $r$ edges,
and all edge and vertex groups trivial.
We obtain an exact sequence
\[
0\to\Z\pi\to(\Z\pi)^r\to{H^1(\pi;\Z\pi)}.
\]
This sequence is in fact the dual of a short free resolution of the augmentation module $\Z$.
If $H^1(\pi;\Z\pi)$ were projective then this sequence would split,
and so the  original resolution for $\Z$ would also split.
This is only possible if $\pi=1$.
Hence  $H^1(\pi;\Z\pi)$ has projective dimension 1 as a right $\Z\pi$-module.

If $m>0$ then we may assume that $\Gamma$ has $m$ vertices and $m+r-1$ edges,  
with the $i$th vertex group being $G_i$, and all edge groups trivial.
Since the vertex groups have one end,
$H^0(G_v;\Z\pi)=H^1(G_v;\Z\pi)=0$ for all $v\in{V}$,
and since the edge groups $G_e$ are trivial,
$H^0(G_e;\Z\pi)=\Z\pi$ for all $e\in{E}$.
Hence $H^1(\pi;\Z\pi)$ is free of rank $m+r-1$.
\end{proof}

We shall use this lemma to match irreducible summands of boundary components
with indecomposable factors of $\pi$.
The next theorem includes parts (2) and (3) of Theorem \ref{C} of the introduction.

\begin{theorem}
\label{*3dg}
Let $(X,\partial{X})$ be an aspherical  $PD_4$-pair such that $\cd\pi=3$,
$H^2(\pi;\Z\pi)=0$ and $\partial{X}=\sqcup_{i=1}^mY_i$.
Let $\nu_i$ be the image  of $\pi_1Y_i$ in $\pi$, for $i\leq{m}$. 
Then
\begin{enumerate}
\item{} $\pi\cong(*_{j=1}^rG_j)*F(s)$,
where $G_j$ is a duality group of dimension $3$,  for all $j\leq{r}$, and $r>0$;
\item{}at most $r+s-1$ of the subgroups $\nu_i$ are not $PD_3$-groups;
\item{} $(X,\partial{X})$ may be obtained by adding mapping cylinders of 
$\Z\nu_i$-homology equivalences over boundary components $Y_i'$
of an aspherical $PD_4$-pair $(X',Y')$ with $\pi_1$-injective boundary.
\end{enumerate}
\end{theorem}

\begin{proof}
 Since $\pi$ is finitely generated and $\cd\pi<\infty$, 
we have $\pi\cong(*_{j=1}^rG_j)*F(s)$, where the factors $G_j$ have one end.
The one-ended factors $G_j$ are duality groups of dimension 3,
since  $\pi$ is of type $FP$,   $H^2(\pi;\Z\pi)=0$ and $\cd\pi\leq3$, 
and there is at least one such factor ($r>0$) since $\cd\pi=3$. 

The inclusion of $\partial{X}$ into $X$ induces isomorphisms
\[
H^q(\pi;\Z\pi)=H^q(X;\Z\pi)\cong{H^q(\partial{X};\Z\pi)}
\]
for $q\leq2$,  by part (1) of Theorem \ref{cd<n}.
Therefore 
 \[
\oplus_{i=1}^mH^1(Y_i;\Z\pi)=H^1(\partial{X};\Z\pi)\cong{H^1(\pi;\Z\pi)}\cong(\Z\pi)^{r+s-1}
\]
and
\[
H^0(\partial{X};\Z\pi)=H^2(\partial{X};\Z\pi)=0.
\]
Hence $H^1(Y_i;\Z\pi)$ is a projective $\Z\pi$-module and $H^2(Y_i;\Z\pi)=0$,
for all $i\leq{m}$.
Since $\Z\pi$ is a free $\Z\nu_i$-module, it follows that 
$H^1(Y_i;\Z\nu_i)$ is a projective $\Z\nu_i$-module and 
$H^2(Y_i;\Z\nu_i)=0$, for all $i\leq{m}$.
Let $\theta_i:\pi_1Y_i\to\nu_i$ be the the epimorphism induced by $inc_i$
and let $\kappa_i=\ker(\theta_i)$ be its kernel.
Then 
\[
H_1(\kappa_i;\Z)\cong{H_1(Y_i;\Z\nu_i)}\cong{H^2(Y_i;\Z\nu_i)}=0,
\]
by Poincar\'e duality for $Y_i$, and so $\kappa_i$ is perfect.

The five term exact sequence of low degree of cohomology with coefficients in 
a $\Z\nu_i$-module $R$ for  $\pi_1Y_i$  as an extension of $\nu_i$ by $\kappa_i$ is
\[
0\to{H^1(\nu_i;R)}\to{H^1(\pi_1Y_i;R)}\to{H^0(\nu_i;H^1(\kappa_i;R))}\to
\]
\[
\to{H^2(\nu_i;R)}\to{H^2(\pi_1Y_i;R)}.
\]
Since $\kappa_i$ maps trivially to $\nu_i$, 
we see that  $H^1(\kappa_i;R)=Hom(\kappa_i,R)$,
and since $\kappa_i$ is perfect,  $H^1(\kappa_i;R)=0$.
Taking $R=\Z\pi$, we see that $H^1(\nu_i;\Z\pi)\cong{H^1(\pi_1Y_i;\Z\pi)}=H^1(Y_i;\Z\pi)$.
Since these are projective $\Z\pi$-modules, 
each $H^1(\nu_i;\Z\pi)$ is in fact free, by Lemma 21.

We may assume that $Y_i\cong(\#_{k=1}^{r_i}Y_{jk})\#E_i\#\Sigma_i$,
where $Y_{ik}$ is aspherical and the image of $\pi_1Y_{ik}$ in $\nu_i$ is nontrivial,  for $k\leq{r_i}$,
$\pi_1E_i\cong{F(s_i)}$ and $\pi_1\Sigma_i\leq\kappa_i$.
Then $\nu_i\cong(*_{k=1}^{r_i}\nu_{ik})*F(s_i)$,
where $\nu_{ik}$ is the image of $\pi_1Y_{ik}$ in $\nu_i$.
by Lemmas \ref{G*Hontopi} and \ref{free factor}.
Since $H^1(\nu_i;\Z\pi)$ is a free module, 
$r_i\geq1$, for $i\leq{m}$,
and since $\oplus_{i=1}^mH^1(\nu_i;\Z\pi)\cong\oplus_{i=1}^mH^1(Y_i;\Z\pi)\cong
{H^1(\pi;\Z\pi)}$,
we have $\Sigma_{i=1}^m(r_i+s_i-1)=r+s-1$.
Hence at most $r+s-1$ of the  subgroups $\nu_i$ have more than one end,
and none are free groups.

The epimorphism from $\pi_1Y_{ik}$ to $\nu_{ik}$ given by restriction
of $\theta_i$ has perfect kernel, by part (4) of Theorem \ref{cd<n}.
Since $Y_{ik}$ is an aspherical  $PD_3$-complex,
it follows that the chain complex for $Y_{ik}$ with coefficients $\Z\nu_{ik}$ 
is a $\Z\nu_{ik}$-resolution of $\Z$,
as in part (7) of Theorem \ref{cd<n}.
Moreover,  $H^3(Y_{ik};\Z\nu_{ik})\cong\Z$.
Hence $\nu_{ik}$ is a $PD_3$-group and $K_{ik}=K(\nu_{ik},1)$, is a $PD_3$-complex.

If $Y_i$ is orientable then we may orient each of the $PD_3$-complexes $K_{ik}$
so that the map from $Y_{ik}$ to $K_{ik}$ determined by the restriction of $\theta_i$
is a degree-1 map. 
Let $P_i=(\#_{k=1}^{r_i}K_{ik})\#E_i$ be the corresponding connected sum.
If $Y_i$ is not orientable then the homotopy type of
$P_i=(\#_{k=1}^{r_i}K_{ik})\#E_i$  is independent of the choices
of orientations of the orientable summands.
In either case, $H_3(\theta_i;\Z^w)$ is an isomorphism 
and so $\theta_i$ lifts to a degree-1 map $f_i:Y_i\to{P_i}$  \cite[Corollary 3]{He77}.
Let $Y'=\sqcup_{i=1}^mP_i$ and let $X'$ be the mapping cylinder of the map from $Y'$ to $X$ 
corresponding to the inclusions of the subgroups $\nu_i$.
Then $X'$ is aspherical and the inclusion of $Y'$ into $X'$ is  $\pi_1$-injective.

Since  $f_i$ is a degree-1 map and and $\kappa_i$ is perfect,
$f_i$ is a $\Z\nu_i$-homology equivalence.
Let $M_i=MCyl(f_i)$ be the mapping cylinder of $f_i$, for $i\leq{m}$,
and let $\breve{X}=X'\cup_{i=1}^mM_i$.
Then $\partial{X}\subset\breve{X}$ and $(X,\partial{X})\simeq(\breve{X},\partial{X})$.
Hence $(X',Y')$ is  a $PD_4$-pair, by  Lemma \ref{excismv}.
\end{proof}

Connected sums of orientable $PD_3$-complexes are not determined by the summands alone, 
as there are choices of local orientation to be made in forming the connected sums.
However, if there is at least one non-orientable summand then the homotopy type of the sum is independent of such considerations. 
If there is a summand which is a non-orientable 3-manifold then this is easy to see
as we may push discs around an orientation-reversing loop. 
It holds in general by virtue of Turaev's work \cite[Chapter 2]{Hi20}.

\section{undoing connected sums}

In this section we shall show that a $PD_4$-pair $(X,\partial{X})$ with non-empty
boundary can be built up from $D^4,S^3$ and $PD_4$-pairs of groups 
by boundary connected sums  and attaching 1-handles if and only if
$X$ is aspherical and $\partial{X}$ is $\pi_1$-injective.
 
The Kurosh Subgroup Theorem gives little information 
on how free factors of subgroups sit in free products; 
the our next lemma is a partial substitute, which is useful in handling free factors.

\begin{lemma}
\cite[Lemma 13.7]{Hi20}
\label{free factor}
Let $B<G$ be groups such that restriction maps $H^1(G;\Z{G})$ onto
$H^1(B;\Z{G})$.
Then every finitely generated free factor of $B$ is also a free factor of $G$.
\end{lemma}

\begin{proof}
The first cohomology group $H^1(G;M)$ of $G$ with coefficients $M$ is 
the quotient of the group of $M$ valued derivations $\Der(G;M)$ by the principal
derivations $\PR(G;M)$ \cite[Exercise III.1.2]{Brown}.
Since restriction clearly maps $\PR(G;M)$ onto $\PR(B;M)$, for any $M$,
the hypothesis implies that restriction maps $\Der(G;\Z{G})$ onto 
$\Der(B;\Z{G}])$.
If $b\in{B}$ generates a free factor of $B$ then there is a derivation 
$\delta:B\to\Z{B}$ such that $\delta(b)=1$ \cite[Corollary IV.5.3]{DD89}.
This may be viewed as a derivation with values in $\Z{G}]$,
and so is the restriction of a derivation $\delta_G:G\to\Z{G}$.
A second application of \cite[Corollary IV.5.3]{DD89} now shows that $b$
generates a free factor of $G$, since $\delta_G(b)=\delta(b)=1$.
If $F=\langle{b_1,\dots,b_s}\rangle$ is a free factor of $B$ then we
may apply this argument to each element $b_i$ of the given basis.
Hence $F$ is a free factor of $G$.
\end{proof}

In Theorem \ref{APD factors} we shall construct a graph from the 
fundamental data of an aspherical $PD_4$-pair $(X,Y)$ (with $Y$ non-empty),
which becomes a template for its assembly from simpler pieces by 
boundary connected sums and attaching 1-handles.
This theorem completes the proof of Theorem \ref{C} of the introduction.

\begin{theorem}
\label{APD factors}
Let $(X,\partial{X})$ be an aspherical $PD_4$-pair, and let $\pi=\pi_1X$.
If  $\pi_1\inc_i$ is injective for $i\leq{m}$
then $(X,\partial{X})$ may be built up from $(D^4,S^3)$ and $PD_4$-pairs 
by boundary connected sums and attaching $1$-handles.
\end{theorem}

\begin{proof}
If $\pi=1$ then $X\simeq*$ and so $(X,\partial{X})\simeq(D^4,S^3)$, 
since the boundary is $\pi_1$-injective.
Thus we may assume that $\pi\not=1$.
Since $X$ is aspherical, $\pi$ is infinite,  so $H_3(\partial{X};\Z\pi)=0$, 
and  no boundary component is $S^3$.
Since the boundary is $\pi_1$-injective, $H^2(\pi;\Z\pi)=0$, by Corollary \ref{pi-inc inj}.
Since $\pi$ is finitely generated and $\cd\pi<\infty$, 
we have $\pi\cong(*_{j=1}^rG_j)*F(s)$, 
where the factors $G_j$  are duality groups of dimension 3,
by part (1) of Theorem \ref{*3dg}.
The free group case $\pi\cong{F(s)}$ is settled in Corollary \ref{freeboundaryinj}
and so we may assume that $r>0$.

We may assume that $\partial{X}=\sqcup_{i=1}^mY_i$ has $m$ connected components $Y_i$.
Let $\nu_i$ be the image of $\pi_1Y_i$ in $\pi$.
Then $\nu_i\cong\sigma_i*F(s_i)$, 
where $\sigma_i$ has $r_i\geq0$ indecomposable factors which are $PD_3$-groups,
and $s_i\geq0$.
Choose elements $b_{ik}\in\pi$, for $k\leq{s_i}$, 
representing a basis for the free factor.
Do this for each $i\leq{m}$.

Since $H^1(\pi;\Z\pi)\cong{H^1(\partial{X};\Z\pi)}$,
by Theorem 2, 
we have an isomorphism $H^1(\pi;\Z\pi)\cong\oplus_{i=1}^m{H^1(\nu_i;\Z\pi)}$.
Applying Lemma \ref{end module}, we see that since $H^1(\pi;\Z\pi)\cong\Z\pi^{r+s-1}$,
the summands $H^1(\nu_j;\Z\pi)$ are all free, so $\sigma_j$ is nontrivial, for $j\leq{m}$.
Moreover, $r+s-1=\Sigma_{i=1}^m(r_i+s_i-1)$.

The isomorphism $H^1(\pi;\Z\pi)\cong\oplus_{i=1}^m{H^1(\nu_i;\Z\pi)}$
implies that $\Der(\pi;\Z\pi)$ maps onto $\oplus_{i=1}^m{\Der(\nu_i;\Z\pi)}$,
as in Lemma \ref{free factor}.
We may choose a derivation $\delta_{ik}\in{\Der(\pi;\Z\pi)}$ which
maps to a derivation $\bar\delta_{ik}\in{\Der(\nu_i;\Z\pi)}$ such that $\bar\delta_{ik}(b_{ik})=1$
and maps trivially to $\Der(\nu_j;\Z\pi)$, for $j\not=i$.
It follows that each of the subgroups $F(s_i)$ is a distinct free factor of $\pi$, 
and so $\Sigma_{i=1}^ms_i\leq{s}$.

Let $\Gamma(\pi,\partial{X})$ be the bipartite graph with vertex set $I\sqcup{J}$ ,
where $I=\{1,\dots,m\}$ and $J=\{\bar1,\dots,\bar{r}\}$,
and an edge from $i\in{I}$ to $\bar{j}\in{J}$ for each indecomposable factor of $\pi_1Y_i$ 
which is a $PD_3$-group and is conjugate to a subgroup of $G_j$.
Then 
\[
\chi(\Gamma(\pi,\partial{X}))=r+m-\Sigma_{i=1}^mr_i=1-s+\Sigma_{i=1}^ms_i\leq1.
\]
Suppose that $\Gamma(\pi,\partial{X})$ is  not connected.
Then there are proper partitions $I=I_<\sqcup{I_>}$ and $J=J_<\sqcup{J_>}$
 such that all edges from $i\in{I_<}$
end in some $j\in{J_<}$ and all edges from $i\in{I_>}$ end in some $j\in{J_>}$.
Let $\pi_<=(*_{j\in{J_<}}G_j)*F(s')$ and let $K_<=K(*_{j\in{J_<}}G_j,1),\vee^{s'}S^1$.
Let $Y_<=\sqcup*_{i\in{I_<}}Y_j$ and let $X_<$ be
the mapping cylinder of the natural map from $Y_<$ to $K_<$. 
Construct another pair $(X_>,Y_>)$ from the remaining terms.
Choose basepoints in $X_<$ and $X_>$ which are not in $Y_<\cup {Y_>}$.
Then the obvious map from $X_<\vee{X_>}$ induces a homotopy equivalence
from $(X_<,Y_<)\vee(X_>,Y_>)$ to $(X,\partial{X})$.
But it is easy to see that a $PD_n$-pair cannot be a non-trivial 1-point union.
(The image of the fundamental class must be 0 on one side,
and so all cap products must be trivial on that side.)
Therefore $\Gamma(\pi,\partial{X})$ is connected.
(If $s=0$ then $s_i=0$, for $i\leq{m}$, so $\chi(\Gamma(\pi,\partial{X}))=1$,
and $\Gamma(\pi,\partial{X})$ is a tree.)

Let $T$ be a maximal tree in $\Gamma(\pi,\partial{X})$.
Let $\lk_T(j)\subset{I}$ be the link of $\bar{j}$ in $T$, for each $j\leq{n}$.
Let $(X_j,Z_j)$ be the $APD_3\partial$-pair with $X_j\simeq{K(G_j,1)}$ 
and $Z_j$ the disjoint union of the indecomposable summands 
of the boundary component $Y_j$  which correspond to edges connecting $i\in{lk_T(j)}$ to $\bar{j}$.
Now form boundary connected sums to obtain an  $APD_3\partial$-pair 
$(\breve{X},\breve{Y})$ with $\pi_1\breve{X}\cong*_{j=1}^rG_j$
and $\breve{Y}=\sqcup_{i=1}^m\breve{Y}_i$,
where $\breve{Y}_i$ is the summand of $Y_i$ with $\pi_1\breve{Y}_i\cong\sigma_i$.
Choose $s_i$ disjoint discs in $Y_i$.
If $w(b_{ik})=1$ we make a boundary connected sum of $(\breve{X},\breve{Y})$ 
with $(D^3\times{S^1},S^2\times{S^1})$ along the $k$th such disc;
if  $w(b_{ik})=-1$ we use $(D^3\tilde\times{S^1},S^2\tilde\times{S^1})$ instead.
Let $(X',Y')$ be the resulting pair, and let $s'=\Sigma_{i=1}^ms_i$.
Then $X'$ is aspherical, $\pi_1X'\cong(*_{j=1}^rG_j)*F(s')$,
and $Y'=\sqcup_{i=1}^mN_i'$, with $N_i'\simeq{N_i}$.
If $s=s'$ then we may realize an isomorphism
$\pi_1X'\cong\pi=\pi_1X$ by a homotopy equivalence of pairs $(X',Y')\simeq(X,\partial{X})$.
In general,
we must also add a 1-handle for each of the $1-\chi(\Gamma(\pi,Y))=s-s'$ edges 
omitted in choosing a maximal tree.
These arise when some $\sigma_i$ has two or more indecomposable factors
which conjugate to subgroups of $G_j$, for the same value of $j$.
We use this fact to guide our choice of 1-handles.

Finally, we invoke Theorem \ref{PD4natural}  and the second assertion of
Lemma \ref{PD4paircriteria} to conclude that the building blocks $(X_j,Z_j)$ are
$PD_4$-pairs of groups.
\end{proof}

\begin{corollary}
Let $G_i$ be a duality group of dimension $3$, for ${j\leq{r}}$.
There is an aspherical  $PD_4$-pair $(X,\partial{X})$ with $\pi_1X\cong(*_{j=1}^rG_j)*F(s)$
and $\pi_1$-injective boundary if and only if $G_j$ is the ambient group of 
a $PD_4$-pair of groups $(G_j,\mathcal{K}_j)$, for $j\leq{r}$.
\qed
\end{corollary}

The case when $\pi$ has no non-trivial free factor is easier;
we shall present it as another  corollary.

\begin{corollary}
\label{pi=*3DG}
Let $(X,\partial{X})$ be a $PD_4$-pair such that $X$ is connected and 
$\partial{X}=\sqcup_{i=1}^mY_i$ is non-empty.
Then $(X,\partial{X})$ is a boundary connected sum of $PD_4$-pairs of groups 
if and only if 
\begin{enumerate}
\item$X$ is aspherical; and
\item$\pi_1X\cong*_{j=1}^nG_j$, where $G_j$ is a duality group of dimension $3$, for $j\leq{n}$;
\item$\pi_1\inc_i$ is a monomorphism, for $i\leq{m}$.
\qed
\end{enumerate}
\end{corollary}

\section{more general groups}

Apart from a few examples,
we know nothing about what might happen when $\cd\pi=3$ and $H^2(\pi;\Z\pi)\not=0$.
We may need to find appropriate generalizations of Lemmas \ref{Kabproj} and  \ref{ddual}.
The argument for Lemma \ref{Kabproj} shows that if $\cd\pi=3$ then 
$K^{ab}$ has projective dimension 1,
but the following example suggests that this may not be useful.

Let $U=T\times[0,1]$.
Then $M=(U\natural{U})\times{S^1}$ is a compact aspherical 4-manifold with 
3 boundary components.
(One of the $\pi_1$-homomorphisms induced by boundary inclusions is an epimorphism,
while the other two are monomorphisms.)
Let $H=\Z^2*\Z^2$ and $J=\Z$.
Then $\pi=\pi_1M\cong{H\times{J}}$ has one end and $\cd\pi=3$.
A K\"unneth Theorem argument gives $H^2(\pi;\Z\pi)\cong{H^1(H;\Z{H})}\otimes{H^1(J;\Z{J})}$.
This is isomorphic to $(\Z{H})^2$ as a right $\Z\pi$-module, 
and so has projective dimension 1.

For a similar example with connected boundary we may replace $U$ by 
the mapping cylinder of the orientation cover of the Klein bottle.

\section{some questions}

Q1. Is every finitely presented group $\pi$ with $\cd\pi=2$ the fundamental group 
of a compact aspherical 4-manifold?
(Can we apply the criterion of Theorem \ref{PD4pair} effectively?)

This might be so even if $\pi$ is a counterexample to the Eilenberg-Ganea conjecture,
i.e., if $\gd\pi=3$. 
(On the other hand, no  counterexample to the Eilenberg-Ganea conjecture
is the fundamental group of a compact aspherical smooth 4-manifold, as observed in \S8.)

\smallskip
Q2.  If $\pi\cong\rho*\sigma$ is the fundamental group of a compact aspherical 4-manifold 
and $\cd\rho=3$ is $\rho$  also such a fundamental group?

\smallskip
Q3.  Which duality groups of dimension $n$ are the ambient groups of $PD_n$-pairs of groups?
This is not known even for $n=3$.

\smallskip
Q4. What can we say if $\pi$ has one end but is not a duality group 
(i.e., $\cd\pi=3$ and $H^2(\pi;\Z\pi)\not=0$)?

\smallskip
Q5.   In \cite{dh1} we conjecture that a $PD_n$-pair $(X,\partial{X})$ with $X$ aspherical 
and $\partial{X}$ a nonempty $(n-1)$-manifold can be realized as a manifold pair $(M,\partial{X})$
(and prove this when $n = 4$ and $\pi_1M$ is elementary amenable).   
What can be said beyond the elementary amenable case?
In \cite{dh1} we give an example of an aspherical $PD_n$-pair which cannot be realized as a manifold pair.   Is there such an example when $n = 4$?


\end{document}